\documentclass[10pt]{article}
\usepackage{amsfonts,latexsym,amstext}
\usepackage{amsmath}
\usepackage{amssymb}
\usepackage[english]{babel}
\usepackage[latin1]{inputenc}

\usepackage[usenames]{color}
\definecolor{red}{rgb}{1.0,0.0,0.0}
\def\red#1{{\textcolor{red}{#1}}}
\definecolor{blu}{rgb}{0.0,0.0,1.0}

\definecolor{gre}{rgb}{0.03,0.50,0.03}

\usepackage[notref,notcite]{showkeys}

\newtheorem{theorem}{Theorem}[section]
\newtheorem{lemma}[theorem]{Lemma}
\newtheorem{proposition}[theorem]{Proposition}
\newtheorem{definition}[theorem]{Definition}

\newtheorem{hypothesis}[theorem]{Hypothesis}

\newtheorem{remark}[theorem]{Remark}
\newtheorem{corollary}[theorem]{Corollary}

\numberwithin{equation}{section}

\newcommand{\myref}[1]{(\ref {#1})}

\def\qed{{\hfill\hbox{\enspace${ \square}$}} \smallskip}
\def\sqr#1#2{{\vcenter{\vbox{\hrule height .#2pt \hbox{\vrule
 width .#2pt height#1pt \kern#1pt \vrule
width .#2pt} \hrule height .#2pt}}}}
\def\square{\mathchoice\sqr54\sqr54\sqr{4.1}3\sqr{3.5}3}

\def\qedo{\hbox{\hskip 6pt\vrule width6pt height7pt
depth1pt  \hskip1pt}\bigskip}

\def\eps{\varepsilon}

\def\ds{\begin{displaystyle}}
\def\eds{\end{displaystyle}}
\def\dis{\displaystyle }
\def\<{\left\langle }
\def\>{\right\rangle }

\def\dim{\noindent \hbox{{\bf Proof.} }}

\def\R{\mathbb R}
\def\N{\mathbb N}

\def\E{\mathbb E}
\def\P{\mathbb P}

\def\calc{{\cal C}}
\def\cald{{\cal D}}

\def\calf{{\cal F}}

\def\calh{{\cal H}}
\def\calk{{\cal K}}
\def\call{{\cal L}}
\def\caln{{\cal N}}
\def\calp{{\cal P}}

\def\calu{{\cal U}}

\def\call{{\cal L}}

\def\to{\rightarrow}

\begin{document}

\title{Stochastic Optimal Control with Delay in the Control I:
\\
solving the HJB equation through partial smoothing }
\date{}

 \author{Fausto Gozzi
\\
Dipartimento di Economia e Finanza,
Universit\`a LUISS - Guido Carli\\
Viale Romania 32,
00197 Roma,
Italy\\
e-mail: fgozzi@luiss.it\\
\\
Federica Masiero\\
Dipartimento di Matematica e Applicazioni, Universit\`a di Milano Bicocca\\
via Cozzi 55, 20125 Milano, Italy\\
e-mail: federica.masiero@unimib.it}

\maketitle
\begin{abstract}
Stochastic optimal control problems governed by delay equations with delay in the control are usually more
difficult to study than the the ones when the delay appears only in the state.
This is particularly true when we look at the associated Hamilton-Jacobi-Bellman (HJB) equation.
Indeed, even in the simplified setting
(introduced first by Vinter and Kwong \cite{VK} for the deterministic case)
the HJB equation is an infinite dimensional second order semilinear Partial Differential Equation (PDE)
that does not satisfy the so-called ``structure condition'' which substantially means
that the control can act on the system modifying its dynamics
at most along the same directions along which the noise acts.
The absence of such condition, together
with the lack of smoothing properties which is a common feature of problems with delay,
prevents the use of the known techniques (based on Backward Stochastic Differential Equations (BSDEs)
or on the smoothing properties of the linear part) to prove the existence of regular solutions of this HJB
equation and so no results on this direction have been proved till now.

In this paper we provide a result on existence of regular solutions of such
kind of HJB equations. This opens the road to prove existence of optimal feedback controls, a task that will be accomplished in the companion paper \cite{FGFM-II}.
The main tool used is a partial smoothing property that
we prove for the transition semigroup associated to the uncontrolled problem.
Such results hold for a specific class of equations and data which arises naturally in many applied problems.
\end{abstract}


\bigskip \noindent


\bigskip \noindent

\textbf{Acknowledgements}:

Financial support from the grant
MIUR-PRIN 2010-11 ``Evolution differential problems: deterministic
and stochastic approaches and their interactions''
is gratefully acknowledged. The second author have been
supported by the Gruppo Nazionale per l'Analisi Matematica,
la Probabilit\`a e le loro Applicazioni (GNAMPA)
of the Istituto Nazionale di Alta Matematica (INdAM).

\newpage

\tableofcontents

\section{Introduction}


Optimal control problems governed by delay equations with delay in the control are usually harder to study than the ones when the delay appears only in the state (see e.g. \cite[Chapter 4]{BDDM07} and \cite{GM,GMSJOTA}). This is true already in the deterministic case but things get worse in the stochastic case.
When one tries to apply the dynamic programming method the main difficulty is the fact that, even in the simplified setting introduced first by Vinter and Kwong \cite{VK} in the deterministic case (see e.g. \cite{GM} for the stochastic case), the associated HJB equation is an infinite dimensional second order semilinear PDE that
does not satisfy the so-called ``structure condition'', which substantially
means that the control can act on the system modifying its dynamics
at most along the same directions along which the noise acts.

The absence of such condition, together with the lack of smoothing properties which is a common feature of problems with delay,
prevents the use of the known techniques, based on BSDE's
(see e.g. \cite{FT2})
or on fixed point theorems in spaces of continuous functions
(see e.g. \cite{CDP1,CDP2,DP3,G1,G2})
or in Gauss-Sobolev spaces
(see e.g. \cite{ChowMenaldi,GGSPA}),
to prove the existence of regular solutions of this HJB equation:
hence no results in this direction have been proved till now.
The viscosity solution technique can still be used (see e.g. \cite{GMSJOTA})
but to prove existence (and possibly uniqueness) of solutions that are merely continuous.
This is an important drawback in this context, since, to prove the existence of optimal feedback control strategies
through the dynamic programming approach, one needs at least the differentiability of the solution in the ``space-like'' variable.

The main aim of this paper is to provide a new result of existence of regular solutions of such HJB equations that holds when the state equation depends linearly on the history of the control and when the cost functional does not depend on such history.
Such results will be exploited in the companion paper \cite{FGFM-II}
to solve the corresponding stochastic optimal control
problem finding optimal feedback control strategies.
This allows to treat satisfactorily a specific class of state equations and data
which arise naturally in many applied problems (see e.g. \cite{PhamBruder,FabbriGozzi08,FedTacSicon,GM,GMSJOTA,Kol-Sha}). 

The key tool to prove such results is the proof of a ``partial'' smoothing property
for the transition semigroup associated to the uncontrolled equation which
we think is interesting in itself and is presented in Section 3.

\medskip
We believe that such tool may allow to treat also examples where the state equation depends
on the history of the state variable, too. To keep things simpler, here we choose to develop and present the result
when this does not happen leaving the extension to a subsequent paper.

\subsection{Plan of the paper}

The plan of the paper is the following:
\begin{itemize}
  \item in Section \ref{section-statement} we give some notations and we present the problem and the main assumptions;
  \item in Section \ref{section-smoothOU} we prove the partial smoothing property for the Ornstein-Uhlenbeck transition semigroup, and we explain how to adapt it to an infinite dimensional setting;
  \item in Section \ref{section-smooth-conv} we introduce some spaces of functions where we will perform the fixed point argument and we prove regularity of some convolutions type integrals;
  \item  in Section \ref{sec-HJB}
we solve the HJB equation
in mild sense.
\end{itemize}

\section{Preliminaries}\label{section-prel}

\subsection{Notation}\label{subsection-notation}

Let $H$ be a Hilbert space. The norm of an element $x$ in $H$ will be
denoted by $\left|  x\right|_{H}$ or simply $\left|  x\right|  $, if no
confusion is possible, and by $\left\langle \cdot,\cdot\right\rangle _H$, or simply
by $\left\langle \cdot,\cdot\right\rangle $ we denote the scalar product in $H$.
We denote by $H^{\ast}$ the dual
space of $H$. Usually we will identify $H$ with its dual
$H^*$. If $K$ is another Hilbert space, $\call(H,K)$ denotes the
space of bounded linear operators from $H$ to $K$ endowed with the usual
operator norm. All Hilbert
spaces are assumed to be real and separable.

\noindent In what follows we will often meet inverses of operators which are not
one-to-one. Let $Q\in \call\left(H,K\right)  $. Then $H_{0}=\ker Q$ is a closed subspace of $H$. Let
$H_{1}$ be the orthogonal complement of $H_{0}$ in $H$: $H_{1}$ is
closed, too. Denote by $Q_{1}$ the restriction of $Q$ to $H_{1}:$ $Q_{1}$ is
one-to-one and $\operatorname{Im}Q_{1}=\operatorname{Im}Q$. For $k\in
\operatorname{Im}Q$, we define $Q^{-1}$ by setting
\[
Q^{-1}\left(  k\right)  :=Q_{1}^{-1}\left(  k\right)  .
\]
The operator $Q^{-1}:\operatorname{Im}Q\rightarrow H$ is called the
pseudoinverse of $Q$. $Q^{-1}$ is linear and closed but in general not
continuous. Note that if $k\in\operatorname{Im}Q$, then
$ Q_{1}^{-1}\left(  k\right)$ is the unique element of
\(\left\{ h  :Q\left(  h\right)  =k\right\}
\)
with minimal norm (see e.g. \cite{Z}, p.209).

In the following, by $(\Omega, \calf, \P)$ we denote a complete probability
space, and by $L^2_\calp(\Omega\times[0,T],H)$
the Hilbert space of all predictable processes
$(Z_t)_{t\in[0,T]}$ with values in $H$, normed by
$\Vert Z\Vert^2 _{L^2_\calp(\Omega\times[0,T],H)}=\E\int_0^T\vert Z_t\vert^2\,dt$.

Next we introduce some spaces of functions.
Let $H$ and $Z$ be real separable Hilbert spaces.
By $B_b(H,Z)$ (respectively $C_b(H,Z)$, $UC_b(H,Z)$) we denote the space of all functions
$f:H\rightarrow Z$ which are Borel measurable and bounded (respectively continuous
and bounded, uniformly continuous and bounded).

Given an interval $I\subseteq \R$ we denote by
$C(I\times H,Z)$ (respectively $C_b(I\times H,Z)$)
the space of all functions $f:I \times H\rightarrow Z$
which are continuous (respectively continuous and bounded).
$C^{0,1}(I\times H,Z)$ is the space of functions
$ f\in C(I\times H)$ such that for all $t\in I$
$f(t,\cdot)$ is Fr\'echet differentiable.
By $UC_{b}^{1,2}(I\times H,Z)$
we denote the linear space of the mappings $f:I\times H \to Z$
which are uniformly continuous and bounded
together with their first time derivative $f_t$ and its first and second space
derivatives $\nabla f,\nabla^2f$.

\noindent If $Z=\R$ we do not write it in all the above spaces.

\subsection{$C$-derivatives}\label{subsection-C-directionalderivatives}

We first recall the definition of $C$-directional derivatives
given in \cite{Mas}, Section 2, and in \cite{FTGgrad}.
Here $H$, $K$, $Z$ are Hilbert spaces.

 \begin{definition}
 \label{df4:Gder}
Let $H$, $K$, $Z$ be real Hilbert spaces.
Let $C:K \rightarrow H$ be a bounded linear operator and let $f:H\rightarrow Z$.
 \begin{itemize}
 \item The $C$-directional
 derivative $\nabla^{C}$ at a point $x\in H$ in the direction $k\in K$ is defined
 as:
 \begin{equation}
 \nabla^{C}f(x;k)=\lim_{s\rightarrow 0}
 \frac{f(x+s Ck)-f(x)}{s},\text{ }s\in\mathbb{R},
 \label{Cderivata}
 \end{equation}
 provided that the limit exists.
 \item We say that a continuous
 function $f$ is $C$-G\^ateaux differentiable at a point $x\in H$
 if $f$ admits the $C$-directional derivative in every direction $k\in K$ and
 there exists a linear operator, called the $C$-G\^ateaux differential, $\nabla^{C}f(x)\in\call(K,Z)$, such that
 $\nabla^{C}f(x;k)=\nabla^{C}f(x)k$ for $x \in H$, $k\in K$.
 The function $f$ is $C$-G\^ateaux differentiable on $H$ if it is
 $C$-G\^ateaux differentiable at every point $x\in H$.

 \item We say that $f$ is $C$-Fr\'echet differentiable
 at a point $x\in H$ if it is $C$-G\^ateaux differentiable and if the limit
 in (\ref{Cderivata}) is uniform for $k$ in the unit ball of $K$. In this case
 we call $\nabla^C f(x)$ the $C$-Fr\'echet derivative (or simply the $C$-derivative) of $f$ at $x$.
 We say that $f$ is $C$-Fr\'echet differentiable on $H$ if it is $C$-Fr\'echet differentiable at every point $x\in H$.
 \end{itemize}
 \end{definition}
Note that, in doing the $C$-derivative, one considers only the directions in $H$ selected in the image of $C$.
%
%
When $Z=\R$ we have $\nabla^C f(x) \in K^*$. Usually we will identify $K$ with its dual
$K^*$ so $\nabla^C f(x)$ will be treated as an element of $K$.

If $f:H\to \R$ is G\^ateaux (Fr\'echet) differentiable on $H$ we have that, given any $C$
as in the definition above, $f$ is $C$-G\^ateaux (Fr\'echet) differentiable on $H$
and
$$
\<\nabla^{C}f(x),k\>_{K}  =\<\nabla f(x),Ck\>_{H}
$$
i.e. the $C$-directional derivative is just the usual directional derivative at a point $x\in H$ in direction
$Ck\in H$. Anyway the $C$-derivative,
as defined above, allows us to deal also with functions that are not G\^ateaux differentiable in every direction.

Now we define suitable spaces of $C$-differentiable functions.

\begin{definition}
\label{df4:Gspaces}
Let $I$ be an interval in $\R$ and let $H$, $K$ and $Z$ be
suitable real Hilbert spaces.
\begin{itemize}
\item
We call $C^{1,C}_{b}(H,Z)$   the space of all functions $f:H\to Z$
which admit continuous and bounded $C$-Fr\'echet derivative. Moreover we call
$C^{0,1,C}_b(I\times H,Z)$ the space of functions $f:I\times H\to Z$
belonging to $C_b(I\times H,Z)$ and such that, for every $t\in I$,
$f(t,\cdot )\in C^{1,C}_b(H,Z)$. When $Z=\R$ we omit it.
\item We call $C^{2,C}_{b}(H,Z)$ the space of all functions $f$ in $C^1_b(H,Z)$
which admit continuous and bounded directional second order derivative $\nabla^C \nabla f$;
by $C^{0,2,C}_b(I\times H,Z)$ we denote the space of functions
$f \in C_b(I\times H,K)$ such that for every $t\in I$,
$f(t,\cdot )\in C^{2,C}_b(H,Z)$. When $Z=\R$ we omit it.
\item
For any $\alpha\in(0,1)$ and $T>0$ (this time $I$ is equal to $[0,T]$) we denote by $C^{0,1,C}_{\alpha}([0,T]\times H)$ the space of functions
$ f\in C_b([0,T]\times H,Z)\cap C^{0,1,C}_b((0,T]\times H)$ such that
the map $(t,x)\mapsto t^{\alpha} \nabla^C f(t,x)$
belongs to $C_b((0,T]\times H,K)$.
The space $C^{0,1,C}_{\alpha}([0,T]\times H)$
is a Banach space when endowed with the norm
\[
 \left\Vert f\right\Vert _{C^{0,1,C}_{\alpha}([0,T]\times H)  }=\sup_{(t,x)\in[0,T]\times H}
\vert f(t,x)\vert+
\sup_{(t,x)\in (0,T]\times H}  t^{\alpha }\left\Vert \nabla^C f(t,x)\right\Vert_{K^*}.
\]
When clear from the context we will write simply
$\left\Vert f\right\Vert _{C^{0,1,C}_{\alpha}}$.
\item For any $\alpha\in(0,1)$ and $T>0$ we denote by $C^{0,2,C}_{\alpha}([0,T]\times H)$
the space of functions
$ f\in C_b([0,T]\times H)\cap C^{0,2,C}((0,T]\times H)$ such that for all $t\in(0,T],\,x\in H$ the map
$ (t,x)\mapsto t^{\alpha} \nabla^C\nabla f(t,x)$ is bounded and continuous as a map
from $(0,T]\times H$ with values in $H\times K$.
The space $C^{0,2,C}_{\alpha}([0,T]\times H)$ turns out to be a Banach space if it is endowed with the norm
\begin{align*}
\Vert  & f \Vert _{C^{0,2,C}_{\alpha}([0,T]\times H)  }\\
  &=\sup_{(t,x)\in[0,T]\times H}
  \vert f(t,x)\vert+
  \sup_{(t,x)\in[0,T]\times H}  \left\Vert \nabla f\left(  t,x\right)  \right\Vert _{H^*}
  +\sup_{(t,x)\in[0,T]\times H}  t^{\alpha}
  \left\Vert \nabla^C\nabla f\left(  t,x\right)  \right\Vert _{H^*\times K^*}.
\end{align*}
\end{itemize}
\end{definition}


%

\section{Setting of the problem and main assumptions}\label{section-statement}

\subsection{State equation}\label{subsection-stateequation}

In a complete probability space $(\Omega, \calf,  \P)$
we consider the following controlled stochastic
differential equation in $\R^n$ with delay in the control:
\begin{equation}
\left\{
\begin{array}
[c]{l}%
dy(t)  =a_0 y(t) dt+b_0 u(t) dt +\int_{-d}^0b_1(\xi)u(t+\xi)d\xi+\sigma dW_t
,\text{ \ \ \ }t\in[0,T] \\
y(0)  =y_0,\\
u(\xi)=u_0(\xi), \quad \xi \in [-d,0),
\end{array}
\right.  \label{eq-contr-rit}
\end{equation}
where we assume the following.

\begin{hypothesis}\label{ipotesibasic}
\begin{itemize}
\item[]
  \item[(i)] $W$ is a standard Brownian motion in $\R^k$, and $(\calf_t)_{t\geq 0}$ is the
augmented filtration generated by $W$;
  \item[(ii)] $a_0\in \call(\R^n;\R^n)$, $\sigma$ is in $\call(\R^k;\R^n)$;
  \item[(iii)] the control strategy $u$ belongs to $\calu$ where
$$\calu:=\left\lbrace z\in L^2_{\calp}(\Omega\times [0,T], \R^m):
u(t)\in U \;a.s.\right\rbrace $$
where $U$ is a closed subset of $\R^n$;
  \item[(iv)] $d>0$ (the maximum delay the control takes to affect the system);
  \item[(v)] $b_0 \in \call(\R^m;\R^n)$;
  \item[(vi)]
$b_1\in L^2([-d,0],\call(\R^m;\R^n))$.
($b_1$ is the density of the time taken by the control to affect the system).
\end{itemize}
\end{hypothesis}

\medskip

Notice that  assumption (vi) on $b_1$ does not cover the
case of pointwise delay since it is technically complicated to
deal with: indeed it gives rise, as we are going to see in next subsection,
to an unbounded control operator $B$, for this reason we leave the extension
of our approach to this case for further research.

\begin{remark}\label{rm:diminf}
Our results can be generalized to the case when the process $y$ is infinite dimensional.
More precisely, let $y$ be the solution of the following controlled stochastic
differential equation in an infinite dimensional Hilbert space $H$, with delay in the control:
\begin{equation}
\left\{
\begin{array}
[c]{l}%
dy(t)  =A_0 y(t) dt+B_0 u(t) dt +\int_{-d}^0 B_1(\xi)u(t+\xi)d\xi+\sigma dW_t
,\text{ \ \ \ }t\in[0,T] \\
y(0)  =y_0,\\
u(\xi)=u_0(\xi), \quad \xi \in [-d,0).
\end{array}
\right.  \label{eq-contr-rit-inf}
\end{equation}
Here $W$ is a cylindrical Wiener process in another Hilbert space $\Xi$, and $(\calf_t)_{t\geq 0}$ is the augmented filtration generated by $W$.
$A_0$ is the generator of a strongly continuous semigroup in $H$.
The diffusion term $\sigma$ is in $\call(\Xi;H)$ and is such that for every $t>0$
the covariance operator
\[
\mathcal{Q}_t^0:=\int_{0}^{t}e^{sA_0}\sigma\sigma^*e^{sA_0^{\ast}}ds
\]
of the stochastic convolution
\[
\int_0^t e^{(t-s)A_0}\sigma\,dW_s
\]
is of trace class and, for some $\gamma \in (0,1)$,
$$
\int_{0}^{t}s^{-\gamma}{\rm Tr}e^{sA_0}\sigma\sigma^*e^{sA_0^{\ast}}ds < + \infty.
$$
The control strategy $u$ belongs to $L^2_{\calp}(\Omega\times [0,T], U_1)$,
where $U_1$ is another Hilbert space, and the space of admissible controls $\calu$ is built in analogy with the finite dimensional case requiring control strategies to take values in a given closed subset $U$ of $U_1$.
On the control operators we assume $B_0 \in \call(U_1;H)$, $B_1: [-d,0] \to\call(U_1,H))$
such that $B_1 u \in L^2([-d,0],H)$ for all $u \in U$.
In this case, following again \cite{GM,VK}, the problem can be reformulated as an abstract evolution equation
in the Hilbert space $\calh$ that this time turns out to be $H\times L^2([-d,0],H)$.
All the results of this paper hold true in this case, under suitable minor changes that will be clarified along the way.
\hfill\qedo
\end{remark}

\subsection{Infinite dimensional reformulation}
\label{subsection-infdimref}

Now, using the approach of \cite{VK} (see \cite{GM} for the stochastic case), we reformulate
equation (\ref{eq-contr-rit}) as an abstract stochastic
differential equation in the Hilbert space $\calh=\R^n\times L^2([-d,0],\R^n)$.
To this end we introduce the operator $A : \cald(A) \subset \calh
\rightarrow \calh$ as follows: for $(y_0,y_1)\in \calh$
\begin{equation}\label{A}
A(y_0 ,y_1 )=( a_0 y_0 +y_1(0), -y_1'), \qquad
\cald(A)=\left\lbrace(y_0,y_1)\in \calh:y_1\in W^{1,2}([-d,0],\R^n), y_1(-d)=0 \right\rbrace.
\end{equation}
We denote by $A^*$ the adjoint operator of $A$:
\begin{equation}
 \label{Astar}
A^{*}_0(y_0 ,y_1 )=( a_0 y_0, y_1'), \qquad
\cald(A^{*})=\left\lbrace (y_0,y_1)\in \calh:y_1\in W^{1,2}([-d,0],\R^n), y_1(0)=y_0 \right\rbrace .
\end{equation}
We denote by $e^{tA}$ the $C_0$-semigroup generated by $A$: for
$y=(y_0,y_1)\in \calh$,
\begin{equation}
e^{tA} \left(\begin{array}{l}y_0 \\y_1\end{array}\right)=
\left(
\begin{array}
[c]{ll}%
e^{ta_0 }y_0+\int_{-d}^{0}1_{[-t,0]} e^{(t+s)a_0 } y_1(s)ds \\[3mm]
y_1(\cdot-t)1_{[-d+t,0]}(\cdot).
\end{array}
\right)  \label{semigroup}
\end{equation}
Similarly, denoting by $e^{tA^*}=(e^{tA})^*$ the $C_0$-semigroup generated by $A^*$,
we have for
$z=\left(z_0,z_1\right)\in \calh $
\begin{equation}
e^{tA^*} \left(\begin{array}{l}z_0 \\z_1\end{array}\right)=
\left(
\begin{array}[c]{ll}
e^{t a_0^* }z_0 \\[3mm]
e^{(\cdot+t) a_0^* }z_0 1_{[-t,0]}(\cdot) +z_1(\cdot+t)1_{[-d,-t)}(\cdot).
\end{array}
\right)  \label{semigroupadjoint}
\end{equation}
The infinite dimensional noise operator is defined as
\begin{equation}
 \label{G}
G:\R^{k}\rightarrow \calh,\qquad Gy=(\sigma y, 0), \; y\in\R^k.
\end{equation}
The control operator $B$
is bounded and defined as
\begin{equation}
 \label{B}
B:\R^{m}\rightarrow \calh,\qquad Bu=(b_0 u, b_1(\cdot)u), \; u\in\R^m
\end{equation}
and its adjoint is
\begin{equation}
 \label{B*}
B^*:\calh^* \rightarrow \R^{m},\qquad B^*(x_0,x_1)=
b^*_0 x_0+\int_{-d}^0 b_1(\xi)^*x_1(\xi)d\xi, \; (x_0,x_1)\in\calh.
\end{equation}
Note that, in the case of pointwise delay the last term of the drift in
the state equation \myref{eq-contr-rit} is $u(t-d)$, hence $b_1(\cdot)$ is
a measure: the Dirac delta $\delta_{-d}$.
Hence in this case $B$ is unbounded as it takes values in
$\R^n \times C^*([-d,0],\R^n)$ (here we denote by $C^*([-d,0],\R^n)$ the dual space of $C([-d,0],\R^n)$).
It will be useful to write the explicit expression of the first component
of the operator $e^{tA}B$ as follows
\begin{equation}\label{eq:etAB}
\left(e^{tA}B\right)_0:\R^m \to \R^n,\qquad    \left(e^{tA}B\right)_0 u=
    e^{ta_0}b_0u+ \int_{-d}^0 1_{[-t,0]}e^{(t+r)a_0}\red{b_1(r)u\,dr}, \quad u \in \R^m.
\end{equation}
%
%

Given any initial datum $(y_0,u_0)\in \calh$ and any admissible control $u\in \calu$ we call $y(t;y_0,u_0,u)$ (or simply $y(t)$ when clear from the context) the unique solution (which comes from standard results on SDE's, see e.g.  \cite{IkedaWatanabe} Chapter 4, Sections 2 and 3)
of (\ref{eq-contr-rit}).

Let us now define the process
$Y=(Y_0,Y_1)\in L^2_\calp(\Omega \times [0,T],\calh)$ as
$$
Y_0(t)=y(t), \qquad Y_1(t)(\xi)=\int_{-d}^\xi u(\zeta+t-\xi)b_1(\zeta)d\zeta,
$$
where $y$ is the solution of equation (\ref{eq-contr-rit}), $u$ is the control process in (\ref{eq-contr-rit}).
By Proposition 2 of \cite{GM},
the process $Y$
is the unique solution of the abstract evolution equation
in $\calh$
\begin{equation}
\left\{
\begin{array}
[c]{l}
dY(t)  =AY(t) dt+Bu(t) dt+GdW_t
,\text{ \ \ \ }t\in[ 0,T] \\
Y(0)  =y=(y_0,y_1),
\end{array}
\right.   \label{eq-astr}
\end{equation}
where $y_0=x_0$ and $y_1(\xi)=\int_{-d}^\xi u_0(\zeta-\xi)b_1(\zeta)d\zeta$.
Note that we have $y_1\in L^2([-d,0];\R^n)$\footnote{This can be seen, e.g.,
by a simple application of Jensen inequality and Fubini theorem. }.
Taking the integral (or mild) form of (\ref{eq-astr}) we have
\begin{equation}
Y(t)  =e^{tA}y+\int_0^te^{(t-s)A}B u(s) ds +\int_0^te^{(t-s)A}GdW_s
,\text{ \ \ \ }t\in[ 0,T]. \\
  \label{eq-astr-mild}%
\end{equation}

\subsection{Optimal Control problem}


The objective is to minimize, over all controls in $\calu$,
the following finite horizon cost:
 \begin{equation}\label{costoconcreto}
J(t,x,u)=\E \int_t^T \left(\red{\bar\ell_0(s)}+\bar\ell_1(u(s))\right)\;ds +\E  \bar\phi(x(T)).
\end{equation}
where \red{$\bar\ell_0:[0,T]\rightarrow \R$} and
$\bar\phi:\R^n\rightarrow \R$ are continuous and bounded while
$\bar\ell_1:U\rightarrow\R$ is measurable and bounded from below.
Referring to the abstract formulation (\ref{eq-astr}) the cost in (\ref{costoconcreto}) can be rewritten also as
\begin{equation}\label{costoconcreto1}
J(t,x;u)=\E \left(\int_t^T \left[\red{\ell_0(s)}+\ell_1(u(s))\right]\,ds + \phi(Y(T))\right),
\end{equation}
where
$\red{\ell_0:[0,T]\rightarrow \R}$, $\ell_1:U\to \R$ are defined by setting
\begin{equation}\label{l_0}
\red{\ell_0(t):=\bar\ell_0(t)} ,
\end{equation}
\begin{equation}\label{l_1}
\ell_1:=\bar\ell_1
\end{equation}
(here we cut the bar only to keep the notation homogeneous)
while $\phi :\calh\rightarrow \R$ is defined as
\begin{equation}\label{fi0}
\phi(x):=\bar\phi(x_0) \quad
\forall x=(x_0,x_1)\in \calh.
\end{equation}

Clearly, under the assumption above,
$\ell_0$ and $\phi$
are continuous and bounded while $\ell_1$ is measurable and bounded from below.
The value function of the problem is
\begin{equation}\label{valuefunction}
 V(t,x):= \inf_{u \in \calu}J(t,x;u).
\end{equation}
We define the Hamiltonian in a modified way, indeed, for $p\in \calh$, $u \in U$,
we define the current value Hamiltonian $H_{CV}$ as
\begin{equation*}
 H_{CV}(p\,;u):=\<p,u\>_{\R^m}+\ell_1(u)
\end{equation*}
%
and the (minimum value) Hamiltonian by
\begin{equation}\label{psi1}
H_{min}(p)=\inf_{u\in U}H_{CV}(p\,;u),
 \end{equation}
The associated HJB equation with unknown $v$ is then formally written as
\begin{equation}\label{HJBformale1}
  \left\{\begin{array}{l}\dis
-\frac{\partial v(t,x)}{\partial t}=\frac{1}{2}Tr \;GG^*\nabla^2v(t,x)
+ \< Ax,\nabla v(t,x)\>_\calh +\red{\ell_0(t)}+ {H}_{min} (\nabla^Bv(t,x)),\qquad t\in [0,T],\,
x\in D(A),\\
\\
\dis v(T,x)=\phi(x).
\end{array}\right.
\end{equation}

To get existence of mild solutions of \myref{HJBformale1} we will need the following assumption.

\begin{hypothesis}\label{ipotesicostoconcreto}
\begin{itemize}
\item[]
  \item[(i)] $\phi\in C_b(\calh)$ and it is given by \myref{fi0} for
  a suitable $\phi \in  C_b(\R^n)$;
  \item[(ii)] $\red{\ell_0\in C_b([0,T])}$ and it is given by \myref{l_0};
  \item[(iii)] $\ell_1:U\rightarrow\R$ is measurable and bounded from below;
  \item[(iv)] the Hamiltonian $H_{min}:\R^m \to \R$ is Lipschitz continuous so
  there exists $L>0$ such that
  \begin{equation}\label{eq:Hlip}
    \begin{array}{c}
\vert H_{min }(p_1)-H_{min }(p_2)\vert\leq L \vert p_1-p_2\vert
\quad \forall\,p_1,\,p_2\in\R^m;
       \\[1.5mm]
\vert H_{min }(p)\vert\leq L(1 + \vert p\vert )
\quad \forall\,p\in\R^m.
    \end{array}
\end{equation}
\end{itemize}
\end{hypothesis}
To get more regular solutions (well defined second derivative $\nabla^B\nabla$,
which will be used to prove existence of optimal feedback controls) we will need
the following further assumption.

\begin{hypothesis}\label{ipotesicostoconcretobis} The Hamiltonian $H_{min}:\R^m \to \R$ is continuously differentiable
  and, for a given $L>0$, we have, beyond \myref{eq:Hlip},
  \begin{equation}\label{eq:Hlipder}
    \begin{array}{c}
\vert \nabla H_{min }(p_1)-\nabla H_{min }(p_2)\vert\leq L \vert p_1-p_2\vert
\quad \forall\,p_1,\,p_2\in\R^m;
    \end{array}
\end{equation}
\end{hypothesis}

\begin{remark}
The assumption \myref{eq:Hlip} of Lipschitz continuity of $H_{min }$ is satisfied e.g. if the set $U$ is compact.
Indeed, for every $p_1,\,p_2\in\R^m$
\[
 \vert H_{min}(p_1)-H_{min}(p_2)\vert \leq \vert \<p_1,u\>-\<p_2,u\>\vert,\quad u\in U
\]
and in the case of $U$ compact the Lipschitz property immediately follows.
The Lipschitz continuity of $H_{min }$ is satisfied also in the case
when $U$ is unbounded, if the current cost has linear growth at infinity.

Moreover the assumption \myref{eq:Hlipder} of Lipschitz continuity of $\nabla H_{min }$ is verified e.g. if the function $\ell_1$ is convex, differentiable with invertible derivative and with $(\ell_1')^{-1}$ Lipschitz continuous since in this case
$(\ell_1')^{-1}(p)=\nabla H_{min}(p)$.
\hfill\qedo
\end{remark}

\begin{remark}\label{remark:costoconcreto}
We list here, in order of increasing difficulty,
some possible generalization of the above assumptions and of the consequent results.
\begin{itemize}
  \item[(i)]\red{
All our results on the HJB equation and on the control problem
could be extended without difficulties to the case
when the boundedness assumption on $\bar \phi$ 
(and consequently on $\phi$) can be replaced by a polynomial growth assumption:
namely that, for some $N\in \N$, the function
\begin{equation}\label{eq:polgrowthphil0}
 x\mapsto \dfrac{\phi(x)}{1+\vert x\vert ^N},
\end{equation}
is bounded.}
The generalization of Theorem \ref{lemma-der-gen} to this case can be achieved by straightforward changes in the proof, on the line of what is done, in a different context, in
\cite{Ce95} or in \cite{Mas-inf-or}.
  \item[(ii)]
Since our results on the HJB equation are based on smoothing properties
(proved in Section \ref{section-smoothOU}) which holds also for
measurable functions, we could consider current cost and final cost only measurable
instead of continuous. The proofs would be very similar but using different
underlying spaces.
  \item[(iii)]
Using the approach of \cite{G2} it seems possible to relax the Lipschitz assumptions on the Hamiltonian function asking only local Lipschitz continuity of the Hamiltonian function,
but paying the price of requiring differentiability of the data.
\end{itemize}
In this paper we do not perform all such generalizations
since we want to concentrate on the main point: {\em the possibility of solving
the HJB equation and the control problem without
requiring the so-called {\bf structure condition}.}
\end{remark}

\section{Partial smoothing for the Ornstein-Uhlenbeck
semigroup}
\label{section-smoothOU}
This section is devoted to what we call the ``partial'' smoothing property
of Ornstein-Uhlenbeck transition semigroups.
First, in Subsection \ref{section-smoothOU-generalsetting}, we give two
results (Theorem \ref{lemma-der-gen} for the first $C$-derivative
and Proposition \ref{lemma-reg-R_t} for the second derivative)
for a general Ornstein-Uhlenbeck transition semigroup in a real separable Hilbert space $H$.
Then in Subsection \ref{section-smoothOU-particular}, we prove two specific results for our problem (Propositions \ref{cor-der} and \ref{lemmaderhpdeb}).

\subsection{Partial smoothing in a general setting}
\label{section-smoothOU-generalsetting}
Let $H, \Xi$ be two real and separable Hilbert spaces and let us consider
the Ornstein-Uhlenbeck process $X^x(\cdot)$ in $H$ which solves the
following SDE in $H$:
\begin{equation}
\left\{
\begin{array}
[c]{l}%
dX(t)  =AX(t) dt+GdW_t
,\text{ \ \ \ }t\ge 0\\
X(0)  =x,
\end{array}
\right.\label{ornstein-gen}
\end{equation}
where $A$ is the generator of a strongly continuous semigroup in
$H$, $(W_t)_{t\geq 0}$ is a cylindrical Wiener process in $\Xi$ and
$G:\Xi\rightarrow H$.
In mild form, the Ornstein-Uhlenbeck process $X^x$ is given by
\begin{equation}
X^x(t)  =e^{tA}x +\int_0^te^{(t-s)A}GdW_s
,\text{ \ \ \ }t\ge 0. \\
  \label{ornstein-mild-gen}
\end{equation}
$X$ is a Gaussian process, namely for every $t>0$, the law of
$X(t)$ is $\caln (e^{tA}x,Q_t)$, the Gaussian measure with mean $e^{tA}x$ and
covariance operator $Q_t$,
where
\[
 Q_t=\int_0^t e^{sA}GG^*e^{sA^*}ds.
\]
The associated Ornstein-Uhlenbeck transition semigroup $R_t$, is defined by setting,
for every $f\in B_b(H)$ and $x\in H$,
\begin{equation}
 \label{ornstein-sem-gen}
R_t[f](x)=\E f(X^x(t))
=\int_K f(z+e^{tA}x)\caln(0,Q_t)(dz).
\end{equation}
where by $X^x$ we denote the Ornstein-Uhlenbeck process above with initial datum
given by $x\in H$.

It is well known (see e.g. \cite[Section 9.4]{DP1} or \cite[Definition 1.159]{FabbriGozziSwiech}),
that $R_t$ has the strong Feller property (i.e. it transforms
bounded measurable functions in continuous ones) at $t>0$ if and only if
\begin{equation}\label{eq:inclusionDZ}
\operatorname{Im}e^{tA}\subseteq \operatorname{Im} Q_t^{1/2},
\end{equation}
and that such property is equivalent to the so-called null-controllability at $t$
of the linear control system identified by the couple of operators
$(A,G)$ (here $z(\cdot)$ is the state and $a(\cdot)$
is the control):
$$
z'(t)=Az(t)+Ga(t),\qquad z(0)=x.
$$
(see e.g. \cite[Appendix B]{DP1}).
Under \myref{eq:inclusionDZ} $R_t$ also transforms any bounded measurable
function $f$ into a Fr\'echet differentiable one,
the so-called ``smoothing'' property,
and
$$
\|\nabla R_t[f]\|_\infty \le \|\Gamma(t)\|_{\call(H)}\|f\|_\infty
$$
where $\Gamma(t):=Q_t^{-1/2}e^{tA}$.

\red{We observe here that \myref{eq:inclusionDZ} (and hence the strong Feller property and the smoothing property mentioned above) cannot be satisfied for all $t>0$ for the problem we study since the image of $Q_t$ (see Lemma \ref{lemmaQ_t} below) is included in
$\R^n\times{0}$ while, unless $t>r$, the image of $e^{tA}$, in the second component, is strictly bigger than the trivial subspace $\{0\}$.\footnote{\red{We note that, also for the Ornstein-Uhlenbeck transition semigroups related to systems with delay only in the state, the strong Feller property and the smoothing property can be verified, under suitable rank conditions, only when $t$ is greater than the maximum delay $r$ (see e.g. \cite[Chapter 10]{DP2}).}}
However here we are interested to the smoothing property or to a similar regularizing property when $t$ is ``small'', namely $t\rightarrow 0$, since this allows, under suitable
assumptions, to solve the HJB equation associated to our control problem using the approach of mild solutions, see Section \ref{sec-HJB} and, e.g. \cite[Chapter 4]{FabbriGozziSwiech} for a survey on such approach.}


\red{With the above motivation in mind we now come back to our general Ornstein-Uhlenbeck
semigroup given in \myref{ornstein-sem-gen} and look at its directional regularizing properties
for special families of functions. We take} another Hilbert space $K$, a bounded operator $C:K\to H$ and extend the smoothing property in two directions:
searching for $C$-derivatives and
applying $R_t$ to a specific class of bounded measurable functions
(see \cite{Lu} for results in this direction in finite dimension).


Let $P:H\rightarrow H$ be a bounded linear operator and endow $\operatorname{Im}(P)$ with the topology inherited
from $H$;
given any $\bar\phi : \operatorname{Im}(P)\to \R$
\red{Borel} measurable and bounded we define a function $\phi \in B_b(H)$, by setting
\begin{equation}\label{fi-gen}
\phi(x)=\bar\phi(Px) \quad
\forall x\in H.
\end{equation}
We prove that, under further assumptions on the operators $A$, $G$, $C$ and $P$,
the semigroup $R_t$ maps functions $\phi$, defined as in (\ref{fi-gen}),
into $C$-Fr\'echet differentiable functions.
%
%

\begin{theorem}\label{lemma-der-gen} Let $A$ be the generator of a strongly continuous semigroup in $H$.
Let $G:\Xi\rightarrow H$.
Let $K$ be another real and separable Hilbert space and let
$C:K\rightarrow H$ be a linear bounded operator. Let $P\in \call(H)$ and let
$\bar\phi:\operatorname{Im}(P)\rightarrow\R$ be Borel measurable and bounded
and define $\phi :H\rightarrow \R$ as in (\ref{fi-gen}).
Fix $t>0$. Then,
$R_t[\phi]$ is $C$-Fr\'echet differentiable if
\begin{equation}\label{eq:inclusion-iff-P}
 \operatorname{Im} \left(P e^{tA}C\right)\subseteq
 \operatorname{Im} \red{(PQ_t P^*)^{1/2}}.
\end{equation}
In this case we have, for every $k\in K$,
\begin{align}
\langle\nabla^C(R_{t}\left[\phi\right]  )(x)  ,k \rangle_K =
\int_{H}
\bar\phi\left(Pz+Pe^{ tA}x\right)
\<Q_{t}^{-1/2}P e^{tA}C k, Q_t^{-1/2}z\>_H\caln(0,Q_t)(dz).
\label{eq:formulaDRT1-gen-P}
\end{align}
Moreover for every $k\in K$
we have the estimate
\begin{equation*}
\vert \<\nabla^C(R_{t}\left[\phi\right]  )(x)  k\>_K\vert
\leq \Vert \bar\phi\Vert_\infty
\left\Vert Q_{t}^{-1/2}P e^{tA} C
\right\Vert_{\call(K;H)} \;\vert k\vert_K,
\end{equation*}
\end{theorem}


\dim
\color{red}
Let $k\in K$. We compute the incremental ratio
in the direction $Ck$.
\begin{align*}
&
\frac{1}{\alpha}
\left[R_{t}[\phi](x+\alpha Ck)-R_{t}[\phi](x)\right]
\\[2mm]
&  =
\frac{1}{\alpha}
\left[\int_{H}\bar\phi\left(Pz+Pe^{tA}(x+\alpha Ck)\right)
\caln\left(0,Q_{t}\right)(dz)
-\int_{H}\bar\phi\left(Pz+Pe^{tA}x\right)\caln\left(0,Q_{t}\right)(dz)\right]
\\[2mm]
&  =
\frac{1}{\alpha}
\left[\int_{H}\bar\phi\left(Pz+Pe^{tA}x\right)
\caln\left(e^{tA}\alpha Ck ,Q_{t}\right)(dz)
-\int_{H}\bar\phi\left(Pz+Pe^{tA}x\right)\caln\left(0,Q_{t}\right)(dz)\right].
\\[2mm]
&  =
\frac{1}{\alpha}
\left[\int_H\bar\phi\left(z_1+Pe^{tA}x\right)
\caln\left(Pe^{tA}\alpha Ck ,PQ_{t}P^*\right)(dz_1)
-\int_H\bar\phi\left(z_1+Pe^{tA}x\right)\caln\left(0,PQ_{t}P^*\right)(dz_1)\right].
\end{align*}
where in the last passage we have performed the change of variable $Pzz_1$.
\color{black}
By the Cameron-Martin theorem, see
e.g. \cite{DP3}, Theorem 1.3.6, the Gaussian measures
\red{$\caln\left(Pe^{tA}\alpha Ck,PQ_tP^*\right)$ and
$\mathcal{N}\left(0,PQ_{t}P^*\right)$ are equivalent if and only if
$Pe^{tA}\alpha Ck\in\red{\operatorname{Im}(PQ_tP^*)^{1/2}}$}.
In such case,
setting, for $\red{y \in \operatorname{Im}(PQ_tP^*)^{1/2}}$,
\begin{align}
&\red{d(t,y,z_1)     =\frac{d\caln\left(y,PQ_{t}P^*\right)}
{d\mathcal{N}\left(0,PQ_{t}P^*\right)  }(z_1)
  =\exp\left\{  \left\langle (PQ_{t}P^*)^{-1/2}
y,(PQ_{t}P^*)^{-1/2}z_1\right\rangle_{H}
-\frac{1}{2}\left|(PQ_{t}P^*)^{-1/2}y\right|_{H}^{2}\right\}  },
\label{eq:density1}
\end{align}
we have, arguing exactly as in \cite{DP1}, proof of Theorem 9.26,
\red{\begin{align*}
&  \nabla^C(R_{t}[\phi])(x)k=
\lim_{\alpha\rightarrow 0}\frac{1}{\alpha}
\int_{H}\bar\phi\left(Pz+Pe^{tA}x\right)\frac{d(t,Pe^{tA}\alpha Ck,z_1)-1}{\alpha}
\caln(0,(PQ_tP^*))(dz_1)
\\[2mm]
&  =\int_{H}\bar\phi\left(Pz+Pe^{tA}x\right)
\<(PQ_tP^*)^{-1/2}P e^{tA}C k, (PQ_tP^*)^{-1/2}z_1\>_{H}\caln(0,(PQ_tP^*))(dz) 
\end{align*}
which gives (\ref{eq:formulaDRT1-gen-P}).
Consequently
\begin{align*}
 \vert \<\nabla^C(R_{t}\left[\phi\right]  )(x) , k\>_K\vert
&\leq
\Vert \bar\phi\Vert_\infty
\left(\int_{H}\<Q_{t}^{-1/2}P e^{tA} Ck,
(PQ_{t}P^*)^{-1/2}z_1\>^2\caln(0,PQ_tP^*)(dz_1)\right)^{1/2}\\[3mm]\nonumber
&= \Vert \bar\phi\Vert_\infty\Vert( PQ_{t}P^*)^{-1/2} P e^{tA}C \Vert_{\call(K;H}\vert k\vert_K. \nonumber
\end{align*}}
This gives the claim.
\qed

\begin{remark}\label{remark-gen0}
In \cite{DP1}, Remark 9.29, it is showed that the analogous of condition \myref{eq:inclusion-iff-P} is also a necessary condition for the Fr\'echet differentiability of $R_t[\phi]$ for any bounded Borel $\phi$. In this case this is not obvious as we deal with a special class of $\phi$ and so the counterexample provided in such Remark may not belong to this class.
\end{remark}

\begin{remark}\label{remark-gen1}
We consider two special cases of the previous Theorem \ref{lemma-der-gen} that will be useful in next section.
\begin{itemize}
  \item[(i)]
Let $K=H$ and $C=I$. In this case Theorem \ref{lemma-der-gen} gives Fr\'echet differentiability:
for $t>0$ $R_t[\phi]$ is Fr\'echet differentiable if \red{
\begin{equation}\label{eq:inclusion-iff-tutta}
 \operatorname{Im}\left( P e^{tA}\right)\in  \operatorname{Im} (PQ_tP^*)^{1/2}
\end{equation}
 and we have,
for every $h\in H$,
\begin{align}
\langle\nabla(R_{t}[\phi])(x),h \rangle_H =
\int_{H}\bar\phi\left(Pz+Pe^{ tA}x\right)
\<Q_{t}^{-1/2}Pe^{tA} h, Q_t^{-1/2}z_1\>_{H}\caln(0,Q_t)(dz_1).
\label{eq:formulaDRT1-gen-tutta}
\end{align}
Moreover for every $h\in H$
we have the estimate
\begin{equation*}
\vert\< \nabla(R_{t}[\phi])(x),h\>_H\vert
\leq \Vert \bar\phi\Vert_\infty\left\Vert Q_{t}^{-1/2}Pe^{tA}
\right\Vert_{\call(H;H)} \;\vert h\vert_H.
\end{equation*}}
\item[(ii)]
Let $K_0$ and $K_1$ be two real and separable Hilbert spaces and let
$K=K_0\times K_1$ be the product space.
Now, given any $\bar\phi\in B_b(K_0)$, we define, in the same way as in (\ref{fi0}),
a function $\phi \in B_b(K)$, by setting
\begin{equation}\label{fi-genbis}
\phi(k)=\bar\phi(k_0) \quad
\forall k=\left(k_0, k_1 \right)\in K.
\end{equation}
Let $P:K\rightarrow \red{K}$ be the projection on the first component of $K$: for every $k=\left(k_0, k_1 \right)\in K$, $Pk=\red{(k_0,0)}$.
Theorem \ref{lemma-der-gen} says that
$R_t[\phi]$ is $C$-Fr\'echet differentiable for every $t>0$ if \red{
\begin{equation}\label{eq:inclusion-iff}
 \operatorname{Im}(Pe^{tA}C)= \operatorname{Im}\left(\left( e^{tA}C\right)_0,0\right)\in  \operatorname{Im} (PQ_tP^*)^{1/2}, \quad \forall t>0
\end{equation}
 and we have, for every $k\in K$,
\begin{align}\label{eq:formulaDRT1-gen}
&\langle\nabla^C(R_{t}\left[\phi\right]  )(x)  ,k \rangle_K  \nonumber\\
&=\int_{K}\bar\phi\left(z_0+(e^{ tA}x)_0\right)\<(PQ_{t}P^*)^{-1/2}\left(\left( e^{tA} Ck\right)_0,0\right),(PQ_tP^*)^{-1/2}z\>_{K}\caln(0,PQ_tP^*)(dz). 
\end{align}
Moreover for every $k\in K$ we have the estimate
\begin{equation}\label{eq:stimaDRT1-gen}
\vert \<\nabla^C(R_{t}\left[\phi\right])(x),k\>_K\vert
\leq \Vert \bar\phi\Vert_\infty
\left\Vert(P Q_{t}P^*)^{-1/2}
\left(\left( e^{tA} C\right)_0,0
\right)
\right\Vert_{\call(K;K)} \;\vert k\vert_K.
\end{equation}}
\end{itemize}
\end{remark}

\begin{remark}\label{rm:crescitapoli-gen}
In Theorem \ref{lemma-der-gen} we prove the partial smoothing for functions $\phi$
defined as in \myref{fi-gen} for functions $\bar \phi$
Borel measurable and bounded. The boundedness assumption on $\bar \phi$
(and consequently on $\phi$) can be replaced by a polynomial growth assumption:
namely that, for some $N\in \N$,
\[
 x\mapsto \dfrac{\bar\phi(x)}{1+\vert x\vert ^N}
\]
is bounded.
The generalization of Theorem \ref{lemma-der-gen} to this case can be achieved
by straightforward changes in the proof, on the line of what is done in
\cite{Ce95} or in \cite{Mas-inf-or}.
\end{remark}

\subsubsection{Second derivatives}

We now prove that, if $\phi$ is more regular, also $\nabla^C R_t[\phi]$ and $\nabla R_t[\phi]$ have more regularity.
This fact, in the context of our model (see Subsection \ref{section-smoothOU-particular}),
will be used in Section \ref{sec-HJB} to prove $C^2$ regularity of
the solution of the HJB equation.


\begin{proposition}\label{lemma-reg-R_t} Let $A$ be the generator of a strongly continuous semigroup in $H$.
Let $G:\Xi\rightarrow H$.
Let $K$ be another real and separable Hilbert space and let
$C:K\rightarrow H$ be a linear bounded operator. Let $P\in \call(H)$ and let
$\bar\phi:\operatorname{Im}(P)\rightarrow\R$ be Borel measurable and bounded
and define $\phi :H\rightarrow \R$ as in (\ref{fi-gen}).
Fix $t>0$.
Assume that (\ref{eq:inclusion-iff-P})
holds true.
If $\bar\phi$ is such that $\phi\in C^1_b(H)$, then
for every $t>0$ the first order derivatives $\nabla^{C}R_{t}\left[\phi\right]$ and
$\nabla R_{t}\left[\phi\right]$ exist and are bounded, with the second one
given by \red{
\begin{align}
&\<\nabla(R_{t}\left[\phi\right]  )(x) , h\>_H= \int_H \<e^{tA^*}\nabla \phi(z+e^{tA}x),h\>\caln(0,Q_t)(dz)=
\nonumber\\
&R_{t}\left[\<\nabla\phi, e^{tA}h\>_H\right](x)
= \int_H \<\nabla \bar\phi(Pz+Pe^{tA}x),Pe^{tA}h \>\caln(0,Q_t)(dz),\
 \quad \forall\, h\in H.
\label{eq:derivateprimelisce}
\end{align}
Moreover the second order derivatives $\nabla\nabla^{C}R_{t}\left[\phi\right]$,
$\nabla^{C}\nabla R_{t}\left[\phi\right]$ exist, coincide, and we have
\begin{align}
&\< \nabla \nabla^{C}(R_{t}\left[\phi\right]  )(x)  k,h\>_H=
\nonumber
\\[2mm]
&=\int_{H} \<\nabla\bar\phi\left(z_1+Pe^{tA}x\right),
Pe^{tA}h\>_{H}
\<(PQ_{t}P^*)^{-1/2}\left( Pe^{tA} Ck\right),
(PQ_tP^*)^{-1/2}z\>_{H}\caln(0,PQ_tP^*)(dz_1).
\label{eq:derivateseconde}
\end{align}
Finally for every $k\in K,\,h\in H$
we have the estimate
\begin{equation}\label{eq:stimaDRT2-gen}
\vert\< \nabla\nabla^C(R_{t}\left[\phi\right])(x)k,h\>_H\vert
\leq \Vert \nabla\bar\phi\Vert_\infty
\Vert (PQ_{t}P^*)^{-1/2}P e^{tA} C
\Vert_{\call(K;H)} \;\vert k\vert_K\,\vert h\vert_H.
\end{equation}}
\end{proposition}
\dim
Formula \myref{eq:derivateprimelisce} is well known (see e.g. \cite[Theorem 4.41-(iii)]{FabbriGozziSwiech}, the last equality follows by our construction of $\phi$ by means of $\bar\phi$.
%
\begin{align*}
 \<\nabla R_{t}\left[\phi\right]  (x),h\>_H
&
=\lim_{\alpha\rightarrow 0}\frac{1}{\alpha}
 \left[\int_H\phi\left(z+e^{tA}(x+\alpha h)\right)
 \caln\left(0,Q_{t}\right)  (dz)
 -\int_H
 \phi\left(z+e^{tA}x\right)  \caln\left(  0,Q_{t}\right)  (dz) \right]
 \\[2mm]
 &  =\int_H\lim_{\alpha\rightarrow 0}\frac{1}{\alpha}
 \left[\phi\left(z+e^{tA}(x+\alpha h)\right)-\phi\left(z+e^{tA}x\right)\right]
 \caln\left(0,Q_{t}\right)  (dz)\\[2mm]
 &  =\int_H\<\nabla \phi\left(z+e^{tA}x\right), e^{tA}h\>_H
 \caln\left(0,Q_{t}\right)  (dz)=R_{t}\left[\<\nabla \phi , e^{tA}h\>_H\right ]  (x).
\end{align*}
The boundedness
of $\<\nabla(R_{t}\left[\phi\right]  )(x) , h\>$ easily follows.
\noindent We compute the second order derivatives starting from
$\nabla \nabla^C R_t[\phi]$. Using \myref{eq:formulaDRT1-gen-P}, \red{performing the change of variable $Pz=z_1$, and applying} the Dominated Convergence Theorem
we get, for $h\in H\,,k\in K$,\red{
\begin{align*}
&\lim_{\alpha\rightarrow 0}\dfrac{1}{\alpha}
\left[\<\nabla^{C}(R_{t}\left[\phi\right]  )(x+\alpha h)  k\>_K
-\<\nabla^{C}(R_{t}\left[\phi\right]  )(x),  k\>_K
\right]=\\[2mm]
&=\lim_{\alpha\rightarrow 0}\dfrac{1}{\alpha}
 \int_{H}\left( \bar\phi\left(Pz+P e^{tA}(x+\alpha h)\right) -
 \bar\phi\left(Pz+Pe^{tA}x\right)\right)
 \<(Q_{t})^{-1/2}\left( Pe^{tA} Ck\right), (Q_t)^{-1/2}z\>_H
 \caln(0,Q_t)(dz)\\[2mm]
&=\lim_{\alpha\rightarrow 0}\dfrac{1}{\alpha}
 \int_{H}\left( \bar\phi\left(z_1+P e^{tA}(x+\alpha h)\right) -
 \bar\phi\left(z_1+Pe^{tA}x\right)\right)
 \<(PQ_{t}P^*)^{-1/2}\left( Pe^{tA} Ck\right), (PQ_tP^*)^{-1/2}z_1\>_{H}\\[2mm] 
 &\qquad \qquad \qquad\qquad  \caln(0,PQ_tP^*)(dz_1)\\[2mm] 
 &= \int_{H} \<\nabla\bar\phi\left(z_1+Pe^{tA}x\right),
 Pe^{tA}h\>_H \<(PQ_{t}P^*)^{-1/2}\left( Pe^{tA} Ck\right), (PQ_tP^*)^{-1/2}z\>_{H}\caln(0,PQ_tP^*)(dz_1),
\end{align*}
Similarly, using (\ref{eq:derivateprimelisce}) and \myref{fi-gen},
we get, for $h\in H\,,k\in K$,
\begin{align*}
&\lim_{\alpha\rightarrow 0}\dfrac{1}{\alpha}
\left[\<\nabla(R_{t}\left[\phi\right]  )(x+\alpha Ck) , h\>_H
-\<\nabla(R_{t}\left[\phi\right]  )(x),  h\>_H
\right]=\\[2mm]
 &=\lim_{\alpha\rightarrow 0}\dfrac{1}{\alpha}
 \int_{H}\left(\<\nabla \bar\phi\left(z_1+Pe^{tA}(x+\alpha Ck)\right),Pe^{tA}h\>_H -
\<\nabla\bar\phi\left(z_1+Pe^{tA}x\right),Pe^{tA}h\>_H\right)
 \caln(0,PQ_tP^*)(dz)\\[2mm]
 &= \int_{H} \<\nabla\bar\phi\left(z_1+Pe^{tA}x\right),Pe^{tA}h\>_{H}
 \<(PQ_{t}P^*)^{-1/2}\left( Pe^{tA}C k\right), (PQ_tP^*)^{-1/2}z_1\>_{H}
 \caln(0,PQ_tP^*)(dz_1).
\end{align*}}
The above immediately implies \myref{eq:derivateseconde} and the estimate \myref{eq:stimaDRT2-gen}.
\qed
%

\subsection{Partial smoothing in our model}
\label{section-smoothOU-particular}
In the setting of Section \ref{section-statement}
we assume that Hypothesis \ref{ipotesibasic} holds true.
We take $\calh=\R^n \times L^2(-d,0;\R^n)$, $\Xi=\R^k$,
$(\Omega, \calf,  \P)$ a complete probability space, $W$ a standard
Wiener process in $\Xi$, $A$ and $G$ as in (\ref{A}) and (\ref{G}).
Then, for $x\in \calh$, we take the Ornstein-Uhlenbeck process $X^x(\cdot)$
given by \myref{ornstein-mild-gen}.
The associated Ornstein-Uhlenbeck transition semigroup $R_t$ is defined
as in \myref{ornstein-sem-gen} for all $f\in B_b(\calh)$.


The operator $P$ of the previous subsection here is the projection \red{$P_0$} on the first component of the space $\calh$, similarly to Remark \ref{remark-gen1}- (ii). \red{Namely for every $x=\left( \begin{array}{l}
x_0\\x_1
\end{array}\right)$, $Px=\left( \begin{array}{l}
x_0\\0
\end{array}\right)$. For the sake of convenience we introduce also the operator $\Pi_0: \calh \rightarrow \R^n$, given by $\Pi_0x=x_0$.}
Hence, given any $\bar\phi\in B_b(\R^n)$, we define, as in (\ref{fi0})
a function $\phi \in B_b(\calh)$, by setting
\begin{equation}\label{fi}
\phi(x)=\bar\phi(\Pi_0 x)=\bar\phi(x_0) \quad
\forall x=\left(
x_0 ,x_1  \right)\in \calh.
\end{equation}
For such functions, the Ornstein-Uhlenbeck semigroup $R_t$ is written as
\begin{equation}
 \label{ornstein-sem-spec}
R_t[\phi](x)=\E \phi (X^x(t))=\E \bar\phi ((X^x(t))_0)=\int_\calh\bar\phi((z+e^{tA}x)_0)\caln(0,Q_t)(dz).
\end{equation}
Concerning the covariance operator $Q_t$ we have the following.
\begin{lemma} \label{lemmaQ_t}
Let $A$ be defined in (\ref{A}), let $G$ be defined by (\ref{G})
and let $t \ge 0$.
Let $Q^0_t$ be the selfadjoint operator in $\R^n$ defined as
\begin{equation}\label{eq:Q^0def}
Q^0_t:=\int_0^t e^{sa_0}\sigma\sigma^*
e^{sa_0^*}\, ds.
\end{equation}
Then for every $(x_0,x_1) \in \calh$ we have
\begin{equation}\label{eq:Q^0Q}
Q_t\left(
x_0, x_1
\right) =\left(
Q^0_t x_0, 0
 \right)
\end{equation}
and so
$$
\operatorname{Im} Q_t
=\operatorname{Im} Q_t^0\times \left\lbrace 0\right\rbrace\subseteq \R^n
\times \left\lbrace 0\right\rbrace$$
Hence, for every $\bar\phi\in B_b(\R^n)$ and for the corresponding
$\phi :\calh\rightarrow\R$ defined in (\ref{fi}) we have
\begin{equation}
 \label{eq:ornstein-sem-specbis}
R_t[\phi](x)=\int_{\R^n}\bar\phi \left(z_0+(e^{tA}x)_0 \right)\caln(0,Q_t^0)(dz_0).
\end{equation}
\end{lemma}
\dim Let $\left(x_0,x_1\right)\in \calh$ and $t\ge 0$.
By direct computation we have
\begin{align*}
Q_t\left( \begin{array}{l} x_0 \\x_1
\end{array}
 \right) &=\int_0^t e^{sA}GG^*e^{sA^*}\left( \begin{array}{l}
x_0 \\x_1
\end{array}
  \right) \, ds\\[1.5mm]
 &=\int_0^te^{sA}\left(\begin{array}{ll}
  \sigma\sigma^*
 &0
 \\0&0\end{array}\right)
 e^{sA^*}\left( \begin{array}{l}
 x_0 \\x_1
 \end{array}
  \right) \, ds\\[1.5mm]
 &=\int_0^te^{sA}\left(\begin{array}{l}
  \sigma\sigma^*
 e^{sa_0^*}x_0\\ 0
 \end{array}
  \right) \, ds=\int_0^t\left(\begin{array}{l}e^{s a_0 }\sigma\sigma^*e^{sa_0^*}x_0\\
 0  \end{array}
\right) \, ds
\end{align*}
from which the first claim (\ref{eq:Q^0Q}) follows.
The second claim (\ref{eq:ornstein-sem-specbis}) is immediate.
\hfill\qed
\begin{remark}\label{remarkinfdim1}
The statement of the above lemma holds true (substituting $\R^n$ with the Hilbert space $K_0$
introduced below) also in the following more general setting.
Let $\calh=K_0\times K_1$ where $K_0$ and $K_1$ are both real separable Hilbert spaces.
Let $\Xi$ be another separable Hilbert space (the noise space) and consider
the Ornstein-Uhlenbeck process
\begin{equation}
X(t)  =e^{tA}x +\int_0^te^{(t-s)A}GdW_s
,\text{ \ \ \ }t\ge 0, \\
  \label{ornstein-mildbis}
\end{equation}
where $A$ generates a strongly continuous semigroup on $\calh$, and $G\in \call(\Xi,\calh)$.
Assume that
\begin{itemize}
  \item $G=(\sigma,0)$ for $\sigma \in \call(\Xi,K_0)$ so
$GG^*=\left(\begin{array}{ll}
  \sigma\sigma^*
 &0
 \\0&0\end{array}\right)$ with $\sigma\sigma^* \in \call(K_0)$;
  \item for every $k_0\in K_0$, $t\ge 0$,
\begin{equation}\label{eq:invarianceeta}
e^{tA}(k_0,0)=  \left(e^{tA_0}k_0, 0\right),
    \end{equation}
where $A_0$ generates a strongly continuous semigroup in $K_0$;
\end{itemize}
then the claim still hold.
Indeed in such case we have, for $t\ge 0$, $k_0\in K_0$, $k_1 \in K_1$,
\begin{equation}\label{eq:invarianceetastar}
  \left(e^{tA^*}(k_0,k_1)\right)_0=e^{tA_0^*}k_0,
    \end{equation}
where $A_0^*$ is the adjoint of $A_0$.\footnote{Indeed once we know that
$e^{tA}(k_0,0)_1=0$ then (\ref{eq:invarianceeta}) is equivalent to ask
(\ref{eq:invarianceetastar}).} So, for $t\ge 0$,
$$
Q_t^0k_0=\int_0^t e^{s A_0} \sigma\sigma^*e^{s A_0^*}k_0  ds
$$
and
$$
 Q_t(k_0,k_1)=(Q^0_t k_0 ,0).
$$
This works, in particular, in the case described in Remark \ref{rm:diminf}.
\hfill\qedo
\end{remark}

We now analyze when Theorem \ref{lemma-der-gen} can be applied
in the cases $C=I$ or $C=B$ concentrating on the cases when
the singularity at $t=0^+$ of \red{$\|Q_t^{-1/2}P_0e^{tA}C\|$} is integrable,
as this is needed to solve the HJB equation (\ref{HJBformale}).

\subsubsection{$C=I$}

By Theorem \ref{lemma-der-gen} we have our partial smoothing
(namely (\ref{eq:formulaDRT1-gen}) and (\ref{eq:stimaDRT1-gen}))
for $C=I$ if
\[
\operatorname{Im}\red{P_0} e^{tA}\subseteq \operatorname{Im}(\red{P_0} Q_t \red{P_0^*}) ^{1/2}.
\]
By Lemma \ref{lemmaQ_t} and \myref{semigroup} this implies
\begin{equation}\label{eq:inclusionsmoothingI}
\operatorname{Im}e^{ta_0}\subseteq \operatorname{Im}(Q^0_t)^{1/2}.
\end{equation}
Since, clearly, $e^{ta_0}$ is invertible and
$\operatorname{Im}(Q^0_t)^{1/2}=\operatorname{Im}Q^0_t $,
then \myref{eq:inclusionsmoothingI} is true if and only if the operator
$Q^0_t$ is invertible. On this we have the following result, taken from
\cite{Z}[Theorem 1.2, p.17] and \cite{Seidman}.

%

\begin{lemma}\label{lm:Q0inv}
The operator $Q_t^0$ defined in (\ref{eq:Q^0def})
is invertible for all $t>0$ if and only if
$$
\operatorname{Im} (\sigma,a_0\sigma, \dots , a_0^{n-1}\sigma)=  \R^n.
$$
This happens if and only if
the linear control system identified by the couple $(a_0,\sigma)$
is null controllable.
In this case, for $t \to 0^+$,
$$
\|(Q_t^0)^{-{1/2}}\| \sim   t^{-r-{1/2}}$$
where $r$ is the Kalman exponent, i.e. the minimum $r$ such that
$$
\operatorname{Im} (\sigma,a_0\sigma, \dots , a_0^{r}\sigma)=  \R^n.
$$
Hence $r=0$ if and only if $\sigma$ is onto.
\end{lemma}

We now pass to the smoothing property.

\begin{proposition}\label{cor-der}
Let $A$ and $G$ be defined respectively by (\ref{A}) and (\ref{G}).
Let $\bar\phi:\R^n\rightarrow\R$ be measurable and bounded and define,
as in (\ref{fi}), $\phi :\calh\rightarrow\R$, by setting $\phi(x)=\bar\phi(x_0)$
for every $x=(x_0,x_1)\in \calh$.
Then, if $Q^0_t$ is invertible, we have the following:
\begin{itemize}
  \item [(i)] the function $(t,x) \mapsto R_t[\phi](x)$
  belongs to $C_b((0,+\infty)\times \calh)$. Moreover it is Lipschitz continuous in $x$
  uniformly in $t\in [t_0,t_1]$ for all $0<t_0<t_1<+\infty$.
\item[(ii)] Fix any $t>0$.
$R_t[\phi]$ is Fr\'echet differentiable and we have,
for every $h\in \calh$,
\begin{align}
\langle\nabla(R_{t}\left[\phi\right]  )(x)  ,h \rangle_\calh =
\int_{\R^n}
\bar\phi\left(z_0+(e^{ tA}x)_0\right)
\<(Q^0_{t})^{-1/2}
\left( e^{tA} h\right)_0, (Q^0_t)^{-1/2}z_0\>_{\R^n}\caln(0,Q^0_t)(dz_0).
\label{eq:formulaDRT1}
\end{align}
where $\left( e^{tA} x\right)_0$, $\left( e^{tA} h\right)_0$
are given by (\ref{semigroup}).
Moreover for every $h\in \calh$
we have the estimate
\begin{equation*}
\vert \<\nabla(R_{t}\left[\phi\right]  )(x),  h\>_\calh\vert
\leq \Vert \bar\phi\Vert_\infty
\left\Vert (Q^0_{t})^{-1/2}
\left( e^{tA} \right)_0
\right\Vert_{\call(\calh;\R^n)} \;\vert h\vert_\calh.
 \end{equation*}
Hence for all $T>0$ there exists $C_T$ such that
\begin{equation}\label{stimadertutta}
 \vert  \<\nabla R_{t}\left[\phi\right]  (x),h\>_\calh \vert
 \leq C_T t^{-r -{1/2}} \Vert \bar\phi \Vert_\infty\vert h\vert_\calh,
 \quad t\in [0,T],
 \end{equation}
 where $r$ is the Kalman exponent which is $0$ if
 and only if $\sigma$ is onto.
\item [(iii)] Fix any $t>0$.
$R_t[\phi]$ is $B$-Fr\'echet differentiable and we have, for every $k\in \R^m$,
\begin{align}
&\<\nabla^B(R_{t}\left[\phi\right]  )(x),  k\>_{\R^m}
=\int_{\R^n}
\bar\phi\left(z_0+\left( e^{tA} x\right)_0\right)
\<(Q^0_{t})^{-1/2}
\left( e^{tA} Bk\right)_0, (Q^0_t)^{-1/2}z_0\>_{\R^n}\caln(0,Q^0_t)(dz_0).
\label{eq:formulaDRT1B}
\end{align}
Moreover, for every $k\in \R^m$,
\begin{equation}\label{eq:formulastimaDRT1B}
 \vert \<\nabla^{B}(R_{t}\left[\phi\right]  )(x),k\>_{\R^m}\vert
\leq \Vert \bar\phi\Vert_\infty
\left\Vert (Q^0_{t})^{-1/2}
\left( e^{tA} B\right)_0
\right\Vert_{\call(R^m,\R^n)} \; \vert k\vert_{\R^m}.
\end{equation}
Hence for all $T>0$ there exists $C_T$ such that
 \begin{equation}\label{stimader}
 \vert  \<\nabla^{B}R_{t}\left[\phi\right]  (x),k\>_{\R^m} \vert
 \leq C_T t^{-r-{1/2}} \Vert \bar\phi \Vert_\infty\vert k\vert_{\R^m},
 \quad t\in [0,T],
 \end{equation}
 where $r$ is the Kalman exponent which is $0$ if
 and only if $\sigma$ is onto.
\end{itemize}
\end{proposition}




\dim Point (ii) immediately follows from the
invertibility of $Q_t^0$, the discussion just before Lemma \ref{lm:Q0inv},
and Theorem \ref{lemma-der-gen}.
Point (i) follows from point (ii) and from the continuity of trajectories
of the Ornstein Uhlenbeck process \myref{ornstein-gen} with $A$ and $G$ given by
(\ref{A}) and (\ref{G}).
Point (iii) follows observing that the operator $\Pi_0 e^{tA}B:\R^m\to\R^n$, given in
\myref{eq:etAB} is well defined and hence, thanks to the
invertibility of $Q_t^0$, Theorem \ref{lemma-der-gen} can be applied.
\hfill\qedo


%

\subsubsection{$C=B$}

By Theorem \ref{lemma-der-gen} we have the partial smoothing
\eqref{eq:formulaDRT1-gen} and (\ref{eq:stimaDRT1-gen})
for $C=B$ if\red{
\begin{equation}\label{eq:inclusionsmoothingB}
\operatorname{Im}P_0e^{tA}B\subset \operatorname{Im}Q^{1/2}_t \, \Longleftrightarrow\, 
\operatorname{Im}\Pi_0e^{tA}B\subset \operatorname{Im}(Q^0_t)^{1/2}
\end{equation}}
Since, as proved e.g. in \cite{Z} (Lemma 1.1, p. 18 and formula (2.11) p. 210),
$$
\operatorname{Im}(Q^0_t)^{1/2}=\operatorname{Im}(\sigma,a_0\sigma, \dots a_0^{n-1}\sigma ),
$$
then, using \myref{eq:etAB},
\myref{eq:inclusionsmoothingB} is verified if and only if
\begin{equation}\label{eq:inclusionsmoothingBbis}
\operatorname{Im}\left(e^{ta_0}b_0 +\int_{-d}^0 1_{[-t,0]}e^{(t+r)a_0}\red{b_1(r)\,dr}
\right)\subseteq \operatorname{Im}(\sigma,a_0\sigma, \dots a_0^{n-1}\sigma ).
\end{equation}
We now provide conditions, possibly weaker than the invertibility
of $Q_t$, under which \myref{eq:inclusionsmoothingBbis} is verified and the singularity at $t=0^+$ of $\|Q_t^{-1/2}\Pi_0 e^{tA}B\|$
is integrable.
We first recall the following result (see \cite{Z}, Proposition 2.1, p. 211).
\begin{proposition}\label{propimzab}
If $F_1$ and $F_2$ are linear bounded operators acting between
separable Hilbert spaces $X$, $Z$ and $Y$, $Z$ such that
$\|F_1^* f \| = \| F_2^* f \|$  for any
$f \in Z^*$, then $Im F_1 = Im F_2$ and $\|F_1^{-1}z\| = \|F_2^{-1}z\|$
for all  $z\in Im F_1$.
\end{proposition}


\begin{proposition}\label{lemmaderhpdeb}
Assume that Hypothesis \ref{ipotesibasic} holds. Assume moreover that, either
\begin{equation}\label{eq:hpdebreg}
\operatorname{Im}(e^{ta_0}b_0)\subseteq
\operatorname{Im}\sigma, \; \forall t > 0;
\qquad
\operatorname{Im}b_1(s)\in\operatorname{Im}\sigma,
\quad  a.e.\, \forall s\in[-d,0]
\end{equation}
or
\begin{equation}\label{eq:hpdebregbis}
\operatorname{Im}\left(e^{ta_0}b_0 +\int_{-d}^0 1_{[-t,0]}e^{(t+r)a_0}\red{b_1(r)\,dr}
\right)
\subseteq\operatorname{Im}\sigma,
\quad \forall t>0.
\end{equation}
Then, for any bounded measurable $\phi$ as in (\ref{fi}),
$R_{t}\left[\phi\right]$ is $B$-Fr\'echet differentiable for every
$t>0$, and, for every $h\in \R^m$,
$\<\nabla^{B}(R_{t}\left[\phi\right])(x),k\>_{\R^m}$ is given by \myref{eq:formulaDRT1B}
and satisfies the estimate \myref{eq:formulastimaDRT1B}.
Moreover for all $T>0$ there exists $C_T$ such that
\begin{equation}
\vert\<  \nabla^{B}(R_{t}\left[\phi\right])  (x),k\>_{\R^m}\vert
\leq C_T t^{-{1/2}} \Vert \bar\phi \Vert_\infty \;\vert k\vert_{\R^m}.
\label{eq:stimaderBnew}
\end{equation}
\end{proposition}
\dim
Consider the following linear deterministic controlled system in $\calh$:
\begin{equation}\label{eq:detcontrsyst}
\left\{
\begin{array}{l}
 dX(t)=AX(t)dt+Gu_1(t)dt \\
X(0)=Bh,
\end{array}
\right.
\end{equation}
where the state space is $\calh$, the control space is $U_1=\R^k$, the control strategy is
$u_1\in L^2_{loc}([0,+\infty);U_1)$, the initial point is $Bh$ with $h\in\R^m$.
Define the linear operator
$$
{\mathcal L}^0_t: L^2([0,t];U_1)\to \R^n, \quad  u_1(\cdot)\mapsto
\int_0^te^{a_0 (t-s)}\sigma u_1(s)ds.
$$
Then the first component of the state trajectory is
\begin{equation}\label{sistcontrfin}
X^0(t)=\Pi_0 e^{tA}Bk + {\mathcal L}^0_t u_1
\end{equation}
Hence $X^0$ can be driven to $0$ in time $t$ if and only if
$$
\Pi_0 e^{tA}Bk \in \operatorname{Im}\mathcal L^0_t
$$
In such case, by the definition of pseudoinverse
(see Subsection \ref{subsection-notation}),
we have that the control which brings $X^0$ to $0$ in time $t$ with minimal $L^2$ norm is $(\mathcal L^0_t)^{-1}\Pi_0 e^{tA}Bk$
and the corresponding minimal square norm is
\begin{equation}
\mathcal{E}\left(  t,Bk\right):  =\min\left\{ \int_{0}
^{t}\left\vert u_{1}(s)\right\vert ^{2}ds:X\left(  0\right)
=B k,\text{ }X^0\left(  t\right)  =0\right\}
=\|(\mathcal L^0_t)^{-1}\Pi_0 e^{tA}Bk\|^2_{L^2(0,t;U_1)}. \label{energyB}
\end{equation}
Since for all $z \in \R^n$ we have
$$
\|(Q_t^0)^{1/2}z\|^2_{\R^n}= |\<Q_t^0z,z\>_{\R^n}|=
\|(\mathcal L^0_t)^* z\|^2_{L^2(0,t;U_1)}
$$
then, by Proposition \ref{propimzab}, we get
$$
\operatorname{Im}\left((Q^0_{t})^{{1/2}}\right)=\operatorname{Im}\mathcal L^0_t .
$$
and
$$
\|(\mathcal L^0_t)^{-1}\Pi_0 e^{tA}Bk\|_{L^2(0,t;U_1)}
=\|(Q^0_{t})^{-{1/2}}\Pi_0 e^{tA}Bk\|_{\R^n}.
$$
Hence, by \myref{energyB}, to estimate
$\|(Q^0_{t})^{-{1/2}}\Pi_0 e^{tA}Bk\|_{\R^n}$
it is enough to estimate the minimal energy to steer $X^0$ to $0$ in time $t$.
When (\ref{eq:hpdebreg}) holds we see, by simple computations, that the control
\begin{equation}\label{control}
\bar u_1(s)=-\dfrac{1}{t}\sigma^{-1}e^{sa_0}b_0 k-\sigma^{-1}b_1(-s)k 1_{[-d,0]}(-t),
\quad s \in [0,t],
\end{equation}
where $\sigma^ {-1}$ is the pseudoinverse of $\sigma$,
brings $X^0$ to $0$ in time $t$.
Hence, for a suitable $C>0$ we get
\[
\mathcal{E}\left(  t,Bk\right)\leq
\int_0^t \bar u_1^2(s) ds
\leq C\left(\dfrac{1}{t}+
\Vert b_1 \Vert^2_{L^2 \left([-d,0];\call(\R^m;\R^n)\right)}\right)|k|^2_{\R^m}.
\]
So, for a possibly different constant $C$, we get $\Vert (Q^0_{t})^{-{1/2}}
\left( e^{tA} Bk\right)_0
\Vert_{\call(\R^m;\R^n)}\leq C t^{-\frac{1}{2}}|k|_{\R^n}$ and the estimate is proved.
If we assume (\ref{eq:hpdebregbis}) we can take as a control,
on the line of \cite{Z}, Theorem 2.3-(iii), p.210,
\begin{equation}\label{controlbis}
\hat u_1(s)=
-\sigma^{*}e^{(t-s)a_0^*}(Q_t^0)^{-1}\left(e^{ta_0}b_0 k+
\int_{-d}^0 1_{[-t,0]}e^{(t+r)a_0}b_1(dr)k\right),
\quad s \in [0,t],
\end{equation}
and use that the singularity of the second term as $t\to 0^+$ is still of order $\dfrac{1}{t}$ since (\ref{eq:hpdebregbis}) holds (see e.g. \cite{SY}, Theorem 1).
Once this estimate is proved, the proof of the $B$-Fr\'echet differentiability
is the same as the one of Proposition \ref{cor-der}-(iii).
%
\qed

\begin{remark}\label{remark-derinfdim}
The above results can be generalized to the case, introduced in Remark \ref{rm:diminf} above, when the first component of the space $\calh$ is infinite dimensional. \red{For brevity we denote $\Pi_0(e^{tA})=(e^{tA})_0$, and similarly for denoting the first component of other operators.}
\begin{itemize}
\item For the case $C=I$ the required partial smoothing
holds if we ask, in place of the invertibility of $Q_t^0$,
that, for every $t> 0$,
\begin{equation}\label{eq:hpdiff}
  \operatorname{Im}\left( e^{tA} \right)_0 \subseteq   \operatorname{Im}(Q_t^0)^{{1/2}}
\end{equation}
which would imply that the linear operator $(Q_t^0)^{-{1/2}}\left( e^{tA} \right)_0$
is continuous from $K_1$ into itself. 
\item For the case $C=B$, the required partial smoothing
holds if we ask that, for every $t> 0$, 
\begin{equation}\label{eq:hpBdiff}
  \operatorname{Im}\left( e^{tA} B \right)_0 \subseteq   \operatorname{Im}(Q_t^0)^{{1/2}}
\end{equation}
which would imply that the operator $(Q_t^0)^{-{1/2}}\left( e^{tA} B\right)_0$ is continuous from $U$ to $K_1$.
\end{itemize}
Clearly, in this generalized setting the estimates (\ref{stimadertutta}) and
(\ref{stimader}) do not hold any more and they depend
on the specific operators $A$, $B$, $\sigma$.
%
\hfill\qedo
\end{remark}

\section{Smoothing properties of the convolution}
\label{section-smooth-conv}
By Proposition
\ref{lemmaderhpdeb}
we know that if \myref{eq:hpdebreg} or \myref{eq:hpdebregbis} hold and $\phi$ is as in (\ref{fi}) with $\bar\phi$ measurable and bounded then $\nabla^{B}(R_t \left[\phi\right]  )(x)$ exists and its norm blows up like
$t^{-{1/2}}$ at $0^+$. Moreover if
$\bar\phi\in C_b(\R^n)$, then $R_t[\phi]\in C_b([0,T]\times\calh)$, see e.g. \cite{PriolaTesi} Proposition 6.5.1 (or the discussion at the end of \cite{PriolaStudia}).

We now prove that, given $T>0$,
for any element $f$ of a suitable family of functions
in $C_b([0,T]\times \calh)$,
a similar smoothing property for the convolution integral
$\int_{0}^{t}R_{t-s} [f(s,\cdot)](x) ds$ holds.
This will be a crucial step to prove the existence and uniqueness
of the solution of our HJB equation in next section.

For given $\alpha \in (0,1)$ we define now a space designed for our purposes.

\begin{definition}\label{df:Sigma}
Let $T>0$, $\alpha \in (0,1)$.
A function $g\in C_b([0,T]\times \calh)$
belongs to $\Sigma^1_{T,\alpha}$ if there exists a function
$f\in C^{0,1}_{\alpha}([0,T]\times \R^n)$ such that
$$g(t,x)=f\left(t,(e^{tA}x)_0\right),
\qquad \forall (t,x) \in [0,T]\times \calh.
$$
\end{definition}
If $g\in \Sigma^1_{T,\alpha}$, for any $t\in(0,T]$ the function $g(t,\cdot)$ is
both Fr\'echet differentiable and $B$-Fr\'echet differentiable.
Moreover, for $(t,x)\in [0,T]\times \calh$, $h \in \calh$, $k\in \R^m$,
$$
\<\nabla g(t,x),h\>_{\calh}=\<\nabla f\left(t,(e^{tA}x)_0\right),(e^{tA}h)_0\>_{\R^n},
\quad and \quad
\<\nabla^B g(t,x),k\>_{\R^m}=\<\nabla f\left(t,(e^{tA}x)_0\right),(e^{tA}Bk)_0\>_{\R^n}.
$$
This in particular imply that, for all $k\in \R^m$
\begin{equation}
\label{eq:nablaperSigma}
\<\nabla^B g(t,x),k\>_{\R^m}=\<\nabla g(t,x),Bk\>_{\R^n \times L^2([-d,0];\R^n)},
\end{equation}
which also means $B^*\nabla g = \nabla^B g$.
For later notational use we call $\bar f\in C_b((0,T]\times \R^n;\R^m)$
the function defined by
$$
\<\bar f(t,y),k\>_{\R^m}=
t^\alpha\<\nabla f\left(t,y\right),(e^{tA}Bk)_0\>_{\R^n},
\qquad (t,y)\in (0,T]\times \R^n, \quad k \in \R^m,
$$
which is such that
$$
t^\alpha\nabla^B g(t,x)=\bar f\left(t,(e^{tA}x)_0\right).
$$
We also notice that if $g\in \Sigma^1_{T,\alpha}$, then in order to have $g$ $B$-Fr\'echet differentiable it suffices to require $(e^{tA}B)_0$
bounded and continuous.

When (\ref{eq:hpdebreg}) or
(\ref{eq:hpdebregbis}) hold we know, by Proposition \ref{lemmaderhpdeb},
that the function $g(t,x)=R_t[\phi](x)$
for $\phi$ given by (\ref{fi}) with $\bar\phi$ bounded and continuous,
belongs to $\Sigma^1_{T,1/2}$.

\begin{lemma}\label{lemma:Sigma}
The set $\Sigma^1_{T,\alpha}$ is a closed subspace of $C^{0,1,B}_{\alpha}([0,T]\times \calh)$.
\end{lemma}
\dim
It is clear that $\Sigma^1_{T,\alpha}$ is a vector subspace of
$C^{0,1,B}_{\alpha}([0,T]\times \calh)$. We prove now that it is closed.
Take any sequence $g_n \to g$ in $C^{0,1,B}_{\alpha}([0,T]\times \calh)$.
Then to every $g_n$ we associate the corresponding $f_n$ and $\bar f_n$.
The sequence $\{f_n\}$ is a Cauchy sequence in
$C_b([0,T]\times \R^n)$. Indeed for any $\epsilon >0$ take
$(t_\epsilon,y_\epsilon)$ such that
$$
\sup_{(t,y) \in [0,T]\times \R^n}|f_n(t,y)-f_m(t,y)|
<
\epsilon + |f_n(t_\epsilon,y_\epsilon)-f_m(t_\epsilon,y_\epsilon)|
$$
Then choose $x_\epsilon \in \calh $ such that
$y_\epsilon=(e^{t_\epsilon A}x_\epsilon)_0$
(this can always be done choosing e.g.
$x_\epsilon=(e^{-t_\epsilon a_0}y_\epsilon,0)$).
Hence we get
$$
\sup_{(t,y) \in [0,T]\times \R^n}|f_n(t,y)-f_m(t,y)|<
\epsilon + |g_n(t_\epsilon,x_\epsilon)-g_m(t_\epsilon,x_\epsilon)|
\le
\epsilon + \sup_{(t,x) \in [0,T]\times \calh}
|g_n(t,x)-g_m(t,x)|.
$$
Since $\{g_n\}$ is Cauchy, then $\{f_n\}$ is Cauchy, too. So there exists a function
$f \in C_b([0,T]\times \R^n)$ such that $f_n \to f$ in $C_b([0,T]\times \R^n)$.
This implies that $g(t,x)=f(t,(e^{tA}x)_0)$ on $[0,T]\times \calh$.
With the same argument we get that there exists a function
$\bar f \in C_b((0,T]\times \R^n;\R^m)$ such that $\bar f_n \to \bar f$
in $C_b((0,T]\times \R^n;\R^m)$.
This implies that $t^\alpha \nabla^B g(t,x)=\bar f(t,(e^{tA}x)_0)$
on $(0,T]\times \calh$.
\qed

Next, in analogy to what we have done defining $\Sigma^1_{T,\alpha}$,
we introduce a subspace $\Sigma^2_{T,\alpha}$ of functions
$g \in C^{0,2,B}_{\alpha}([0,T]\times \calh)$
that depends in a special way on the variable $x\in\calh$.
\begin{definition}\label{df:Sigma2}
A function $g\in C_b([0,T]\times \calh)$ belongs to $\Sigma^2_{T,\alpha}$
if there exists a function
$f\in C^{0,2}_{\alpha}([0,T]\times \R^n)$ such that
for all $(t,x) \in [0,T]\times \calh$,
\begin{align*}
&g(t,x)=f\left(t,(e^{tA}x)_0\right).
\end{align*}
\end{definition}
If $g\in \Sigma^2_{T,\alpha}$ then for any $t\in(0,T]$ the function $g(t,\cdot)$
is Fr\'echet differentiable and
$$
\<\nabla g(t,x),h\>_\calh=
\<\nabla f\left(t,(e^{tA}x)_0\right),(e^{tA}h)_0\>_{\R^n},
\quad \hbox{for $(t,x)\in [0,T]\times \calh$, $h \in \calh$.}
$$
Moreover also $\nabla g(t,\cdot)$ is $B$-Fr\'echet differentiable and
$$
\<\nabla^B\left(\nabla g(t,x)h\right)),k\>_{\R^m}=
\<\nabla^2 f\left(t,(e^{tA}x)_0\right)(e^{tA}h)_0,(e^{tA}Bk)_0 \>_{\R^n},
\quad \hbox{for $(t,x)\in [0,T]\times \calh$, $h \in \calh$, $k \in \R^m$.}
$$
We also notice that, since the function $f$ is twice continuously
Fr\'echet differentiable the second order derivatives
$\nabla^B\nabla g$ and $\nabla\nabla^B g$ both exist and coincide:
$$
 \<\nabla^B\<\nabla g(t,x),h\>_{\calh},k\>_{\R^m}
=\<\nabla\<\nabla^B g(t,x),k\>_{\R^m},h\>_{\calh}.
$$
Again for later notational use we call
$\bar f_1\in C_b([0,T]\times \R^n;\R^m)$
the function defined by
$$
\<\bar f_1(t,y),h\>_{\R^m}=
\<\nabla f\left(t,y\right),(e^{tA}Bh)_0\>_{\R^n},
\qquad (t,y)\in [0,T]\times \R^n, \quad h \in \R^m,
$$
which is such that
$$
\nabla^B g(t,x)=\bar f_1\left(t,(e^{tA}x)_0\right).
$$
Similarly we call
$\bar{\bar f}\in C_b\left((0,T]\times \R^n;\call(\calh,\R^m)\right)$
the function defined by
$$
\<\<\bar{\bar f}(t,y),h\>_\calh,k\>_{\R^m}=
t^\alpha\<\nabla^2 f\left(t,y\right)(e^{tA}h)_0,(e^{tA}Bk)_0 \>_{\R^n}
\qquad (t,y)\in [0,T]\times \R^n, \quad h \in \calh, \; k\in \R^m,
$$
which is such that
$$
t^\alpha\nabla^B\nabla g(t,x)=t^\alpha\nabla\nabla^B g(t,x)=
\bar{\bar f}\left(t,(e^{tA}x)_0\right).
$$

\red{\begin{lemma}\label{lemma-nuovo}
Assume that Hypothesis \ref{ipotesibasic} holds and let (\ref{eq:hpdebreg}) or (\ref{eq:hpdebregbis}) hold true. Then for $0\leq s \leq t$
\begin{equation}\label{eq:incl-nuova}
\operatorname{Im}(P_0e^{tA}B)\subseteq\operatorname{Im}\Big(P_0e^{sA}Q_{t-s}e^{sA^*}P_0^*\Big), \text{ equivalently } \operatorname{Im}(\Pi_0e^{tA}B)\subseteq\operatorname{Im}\Big(\Pi_0e^{sA}Q_{t-s}e^{sA^*}\Pi_0^*\Big)
\end{equation}
and 
\begin{equation}\label{eq:stima-nuova}
\Vert (Q^1_{t-s})^{-1/2}
\Pi_0e^{tA}B)\Vert \leq \frac{c}{(t-s)^{1/2}}.
\end{equation}
\end{lemma}
\dim 
Notice that the operators $P_0,\, P_0^*:\calh\rightarrow \calh$ coincide, so for every $x=\Big(\begin{array}{l}x_0\\x_1
\end{array}\Big)\in \calh$ we can write 
$$
P_0e^{sA}Q_{t-s}e^{sA^*}P_0 x=\Big(\begin{array}{l}e^{s a_0}Q^0_{t-s}e^{s a_0^*}x_0\\ 0\end{array}\Big)
$$
and so inclusion \eqref{eq:incl-nuova}, is equivalent to prove that for all $k\in \R^m$
\begin{equation}\label{eq:incl-nuova-1}
e^{t a_0}b_0k+\int_{-d}^{0}1_{[-t,0]}e^{(t+r)a_0}b_1(r)k\,dr
\in \operatorname{Im}(e^{s a_0}Q^0_{t-s}e^{s a_0^*}).
\end{equation}
Now we notice that for $z\in \R^n$
\begin{equation*}
\langle e^{s a_0}Q^0_{t-s}e^{s a_0^*}z,z\rangle_{\R^n}=\langle 
Q^0_{t-s}e^{s a_0^*}z,e^{s a_0^*}z\rangle_{\R^n}=\Vert  (\call^0_{t-s})^*e^{s a_0^*}z\Vert ^2
\end{equation*}
where 
$$
{\mathcal L}^0_{t-s}: L^2([0,t];U_1)\to \R^n, \quad  u_1(\cdot)\mapsto
\int_0^{t-s}e^{a_0 (t-s-r)}\sigma u_1(r)dr.
$$
By Proposition \ref{propimzab}
$$
 \operatorname{Im}(e^{s a_0}Q^0_{t-s}e^{s a_0^*})^{1/2}=\operatorname{Im}(e^{s a_0}{\mathcal L}^0_{t-s})
 $$
 Next we define 
 \begin{equation}
\label{Q1L1}
Q^1_{t-s}=e^{s a_0}Q^0_{t-s}e^{s a_0^*}
\qquad{\mathcal L}^1_{t-s}= e^{s a_0}{\mathcal L}^0_{t-s}.
 \end{equation}
We look for a control $u_1\in U_1$ which brings $x_0$ to $0$ in time $t-s$, that is such that
$$
\call^1_{t-s}u_1+e^{t a_0}b_0k+\int_{-d}^{0}1_{[-t,0]}e^{(t+r)a_0}b_1(r)k\,dr.
$$
It turns out that
\begin{equation*}
u_1(r)=-\dfrac{1}{t-s}\sigma^{-1}e^{ra_0}b_0 k-\sigma^{-1}b_1(-r)k 1_{[-d,0]}(-t),
\quad  r\in [0,t-s],
\end{equation*}
is such a control, and we deduce that 
so \eqref{eq:incl-nuova} and estimate \eqref{eq:stima-nuova} hold true.
\qed
}

We now pass to the announced smoothing result.

\begin{lemma}\label{lemma reg-convpercontr}
Let (\ref{eq:hpdebreg}) or (\ref{eq:hpdebregbis}) hold true.
Let $T>0$ and let $\psi:\R^m\rightarrow\R$ be a continuous
function satisfiying Hypothesis (\ref{ipotesicostoconcreto}), estimates
(\ref{eq:Hlip}).
Then
\begin{itemize}
 \item [i)]
for every $g \in \Sigma^1_{T,{1/2}}$, the function
$\hat g:[0,T]\times \calh\rightarrow \R$ belongs to $\Sigma^1_{T,{1/2}}$ where
\begin{equation}\label{iterata-primag}
\hat g(t,x) =\int_{0}^{t}
R_{t-s} [\psi(\nabla^{B}g(s,\cdot))](x)  ds.
\end{equation}
 Hence, in particular, $\hat g(t,\cdot)$
is $B$-Fr\'echet differentiable for every $t\in (0,T]$ and, for all $x\in \calh$,
\begin{equation}\label{stimaiterata-primag}
\left\vert \nabla ^B(\hat g(t,\cdot))(x) \right\vert_{(\R^m)^*}    \leq
    C\left(t^{{1/2}}+\Vert g \Vert_{C^{0,1,B}_{{1/2}}}\right).
 \end{equation}
 If $\sigma$ is onto, then $\hat g(t,\cdot)$
is Fr\'echet differentiable for every $t\in (0,T]$ and, for all $h\in\calh$, $x\in \calh$,
\begin{equation}\label{stimaiterata-primagbis}
\left\vert\nabla(\hat g(t,\cdot))(x) \right\vert_{\calh^*}
   \leq
    C\left(t^{{1/2}}+\Vert g \Vert_{C^{0,1}_{{1/2}}}\right).
\end{equation}
\item [ii)] Assume moreover that $\psi\in C^1(\R^m)$.
For every $g \in \Sigma^2_{T,{1/2}}$, the function
$\hat g$ defined in (\ref{iterata-primag})
belongs to $\Sigma^2_{T,{1/2}}$. Hence, in particular, the second order derivatives
$\nabla\nabla^B\hat g(t,\cdot)$ and $\nabla^B\nabla\hat g(t,\cdot)$ exist, coincide and
for every $t\in (0,T]$ and, for all $x\in \calh$,
\begin{equation}\label{stimaiterata-secondag}
\left\vert \nabla ^B\nabla(\hat g(t,\cdot))(x) \right\vert_{\calh^*\times(\R^m)^*}    \leq
    C
\Vert g \Vert_{C^{0,2,B}_{{1/2}}}
 \end{equation}
 If $\sigma$ is onto, then $\hat g(t,\cdot)$
is twice Fr\'echet differentiable and for every $t\in(0,T]$, for all $h\in\calh$ and $x\in \calh$,
\begin{equation}\label{stimaiterata-secondagbis}
\left\vert\nabla^2(\hat g(t,\cdot))(x) \right\vert_{\calh^*\times \calh^*}
   \leq
    C
\Vert g \Vert_{C^{0,2}_{{1/2}}}
\end{equation}
\end{itemize}
\end{lemma}
\dim
We start by proving that (\ref{iterata-primag}) is $B$-Fr\'echet differentiable
and we exhibit its $B$-Fr\'echet derivative.
Recalling (\ref{ornstein-sem-gen}) we have
\begin{align*}
  \int_{0}^{t}
R_{t-s} \left[\psi\left(\nabla^{B}(g(s,\cdot))\right)\right](x)  ds
&=\int_{0}^{t}
\int_{\calh} \psi\left(\nabla^{B}(g(s,\cdot))
\left(z+e^{(t-s)A}x \right)\right)\caln(0,Q_{t-s})(dz)
\end{align*}
By the definition of $\Sigma^1_{T,{1/2}}$, we see that
\begin{align}\label{eq:nablaBshiftg}
 &s^{1/2} \nabla^{B}g(s,z+e^{(t-s)A}x)
 =
 \bar f\left(s, (e^{sA}z)_0+(e^{tA}x)_0\right)
 \qquad \forall t\ge s>0, \; \forall x,z \in \calh.
\end{align}
Hence the function $\hat f$ associated to $\hat g$ is
\begin{align*}
\hat f (t,y)&=\int_{0}^{t}
\int_{\calh} \psi\left(s^{-{1/2}}
 \bar f\left(s, (e^{sA}z)_0+y\right)
\right)\caln(0,Q_{t-s})(dz)
\end{align*}
with, by our assumptions on $\psi$,
$$
\Vert\hat f \Vert_\infty \le C\int_{0}^{t}\left( 1 +
s^{-{1/2}} \Vert \bar f \Vert_\infty \right)ds
$$
To compute the $B$-directional derivative we look at the limit
\begin{align*}
 &\lim_{\alpha\rightarrow 0}\dfrac{1}{\alpha} \left[\int_{0}^{t} R_{t-s}\left[\psi\left( \nabla^{B}(g(s,\cdot))\right)\right](x+\alpha Bk)ds -\int_{0}^{t}R_{t-s}\left[\psi\left( \nabla^{B}(g(s,\cdot))\right)\right](x)  ds\right].
\end{align*}
From what is given above we get \red{
\begin{align*}
&\int_{0}^{t}R_{t-s}\left[\psi\left( \nabla^{B}(g(s,\cdot))\right)\right](x+\alpha Bk)ds \\
&=\int_{0}^{t} \int_\calh  \psi\left( s^{-{1/2}}
\bar f\left(s, (e^{sA}z)_0+(e^{tA}(x+\alpha Bh))_0\right)
\right)\caln(0,Q_{t-s})(dz) ds
\\[2mm]
&=\int_{0}^{t} \int_{\R^n}
\psi\left(s^{-{1/2}}
\bar f\left(s, z_1+(e^{tA}(x+\alpha Bk)_0\right)
\right)
\caln\left(\left(0, 0\right),\Pi_0e^{sA}Q_{t-s}e^{sA}\Pi_0^*\right)(dz_1)ds
\end{align*}
where in the last passage we have performed the substitution $z_1=(e^{sA}z)_0=\Pi_0 e^{sA}z$. We let  $Q^1_{t-s}=\Pi_0e^{sA}Q_{t-s}e^{sA^*}\Pi_0^*$ already defined in \eqref{Q1L1}, and we get with $z_2=z_1+(e^{tA}\alpha Bk)_0$
\begin{align*}
&\int_{0}^{t}R_{t-s}\left[\psi\left( \nabla^{B}(g(s,\cdot))\right)\right](x+\alpha Bk)ds \\
&=\int_{0}^{t} \int_{\R^n}\psi\left(s^{-{1/2}}\bar f\left(s, z_2+(e^{tA}x)_0\right)\right)
\caln\left((e^{tA}\alpha Bk)_0,Q^1_{t-s}\right)(dz_2)ds\\[2mm]
&=\int_{0}^{t} \int_{\R^n}\psi\left(s^{-{1/2}}
\bar f\left(s, z_2+(e^{tA}x)_0\right)
\right)
d(t,t-s,(e^{tA}\alpha  Bk)_0,z_2)
\caln\left(0,Q^1_{t-s}\right)(dz_2)ds,
\end{align*}
By Lemma \ref{lemma-nuovo}, inclusion \eqref{eq:incl-nuova}, the Gaussian measures $\mathcal{N}\left(\left(e^{t_1 A}y\right)_0,Q^1_{t_2}\right)  $ and $\mathcal{N}\left(0,Q^1_{t_2}\right)  $
where
\begin{align}
&\frac{d\mathcal{N}\left(\left(e^{t_1 A}y\right)_0,Q^1_{t_2}
\right)  }{d\mathcal{N}\left(  0,Q^1_{t_2}\right) }(z_2)
\nonumber
\\[2mm]
&  =\exp\left\{  \left\langle  (Q^1_{t_2})^{-1/2}
\left(e^{t_1A}y\right)_0,(Q^1_{t_2})^{-{1/2}}z_2\right\rangle_{\R^n}
-\frac{1}{2}\left| ( Q^1_{t_2})^{-{1/2}}\left(
 e^{t_1A}y\right)_0\right|_{\R^n}^{2}\right\}  .
\label{eq:density1bisg}
\end{align}
Hence
\begin{align*}
&\lim_{\alpha\rightarrow 0}\dfrac{1}{\alpha}
 \left[\int_{0}^{t}
 R_{t-s}\left[\psi\left( \nabla^{B}(g(s,\cdot))\right)\right](x+\alpha Bk)ds -
\int_{0}^{t}
 R_{t-s}\left[\psi\left( \nabla^{B}(g(s,\cdot))\right)\right](x)  ds\right]=
 \\[2mm]
&=\lim_{\alpha\rightarrow 0}\dfrac{1}{\alpha}
\int_{0}^{t} \int_{\R^n}
\psi\left(s^{-{1/2}}
\bar f\left(s, z_2+(e^{tA}x)_0\right)
\right)
\frac{d(t,t-s,(e^{tA}\alpha Bk)_0,z_2)-1}{\alpha}
\caln\left(0,Q^1_{t-s}\right)(dz_2)ds
\\[2mm]
&=
\int_{0}^{t} \int_{\R^n}
\psi\left(s^{-{1/2}}
\bar f\left(s, z_2+(e^{tA}x)_0\right)
\right)
\<(Q^1_{t-s})^{-{1/2}}\left(
 e^{tA} Bk\right)_0,( Q^1_{t-s})^{-{1/2}}z_2\>_{\R^n}
\caln\left(0,Q^1_{t-s}\right)(dz_2)ds.
\end{align*}
Since the above limit is uniform for $k$ in the unit sphere, then
we get the required $B$-Fr\'echet differentiability and
\begin{align}
&\<\nabla ^B \left(\int_{0}^{t}
R_{t-s}\left[\psi\left( \nabla^{B}(g(s,\cdot)\right)\right]  ds
\right) (x),k \>_{\R^m}=
\label{eq:derBconvnew}
\\[3mm]\nonumber
&=\int_{0}^{t} \int_{\R^n}
\psi\left(s^{-{1/2}}\bar f\left(s, z_2+(e^{tA}x)_0\right)\right)
\<(Q^1_{t-s})^{-{1/2}}\left(e^{tA} Bk\right)_0,( Q^1_{t-s})^{-{1/2}}z_2\>_{\R^n}
\caln\left(0,Q^1_{t-s}\right)(dz_2)ds.
\end{align}
Finally we prove the estimate (\ref{stimaiterata-primag}). Using the above representation, the Holder inequality and Lemma \ref{lemma-nuovo}, estimate \eqref{eq:stima-nuova}, we have
\begin{align*}
&\left\vert\<\nabla ^B \left(\int_{0}^{t}
R_{t-s}\left[\psi\left( \nabla^{B}(g(s,\cdot))\right)\right]  ds
\right) (x),k \>_{\R^m}\right\vert \le
\\[3mm]
&\leq C\int_{0}^{t}
 \int_{\R^n}  \left(1+\left\vert
s^{-{1/2}} \bar f\left(s, z_2+(e^{tA}x)_0\right)
\right\vert\right)
 \left \vert
\<(Q^1_{t-s})^{-{1/2}}
\left( e^{tA} Bk\right)_0, (Q^1_{t-s})^{-{1/2}}z_2\>_{\R^n}
 \right\vert
\caln(0,Q^1_{t-s})(dz_2) ds\\
&\leq C\int_{0}^{t}\left(1+s^{-{1/2}} \left\Vert g \right\Vert_{C^{0,1,B}_{{1/2}}}
\right)
\left\Vert (Q^1_{t-s})^{-{1/2}}(e^{tA}Bk)_0 \right\Vert_{\call(\R^m;\R^n)} ds \\
 &\leq C\int_{0}^{t}\left(1+s^{-{1/2}}\left\Vert g \right\Vert_{C^{0,1,B}_{{1/2}}}\right)
 (t-s)^{-{1/2}}\vert k\vert_{\R^m} \,ds \leq C\left(t^{{1/2}}
 +\left\Vert g \right\Vert_{C^{0,1,B}_{{1/2}}} \right)
\vert k\vert_{\R^m}.
\end{align*}
Observe that in the last step we have used estimate \eqref{eq:stima-nuova}
which follows from the proof of Lemma \ref{lemma-nuovo}}.
Moreover we have also used that
\begin{equation}
\label{eq:betaeuler}
 \int_0^t(t-s)^{-{1/2}}s^{-{1/2}}ds=\int_0^1(1-x)^{-{1/2}}x^{-{1/2}}dx=\beta \left({1/2},{1/2} \right) ,
\end{equation}
where by $\beta(\cdot,\cdot)$ we mean the Euler beta function.

The Fr\'echet differentiability and the estimate
(\ref{stimaiterata-primagbis})
is proved exactly in the same way using the fact that $\sigma$ is onto
and Proposition \ref{cor-der}.
Now we consider the case of $g\in\Sigma^2_{T,{1/2}}.$
We start by proving that (\ref{iterata-primag}) is Fr\'echet differentiable
and, in order to compute the Fr\'echet derivative, we use \myref{eq:derivateprimelisce}
(which is true for every $\phi \in C^1_b(H)$)
looking at the limit, for $h \in \calh$,
\begin{align}
&\<\nabla\hat g(t,x),h\>_\calh
 =\lim_{\alpha\rightarrow 0}\dfrac{1}{\alpha}
  \left[\int_{0}^{t}
  R_{t-s}\left[\psi\left( \nabla^{B}(g(s,\cdot))\right)\right](x+\alpha h)ds -
 \int_{0}^{t}
  R_{t-s}\left[\psi\left( \nabla^{B}(g(s,\cdot))\right)\right](x)  ds\right]
  \nonumber
  \\
&=\int_{0}^{t} R_{t-s}
\left[\<\nabla\psi\left( \nabla^{B}(g(s,\cdot))\right),
\nabla \nabla^{B}(g(s,\cdot)) e^{(t-s)A}h\>_{\R^m}\right](x)ds
\nonumber
\\
&=\int_{0}^{t} \int_\calh
\<\nabla\psi\left( \nabla^{B}(g(s,z+e^{(t-s)A}x))\right),
\nabla \nabla^{B}(g(s,z+e^{(t-s)A}x)) e^{(t-s)A}h\>_{\R^m}
\caln\left(0,Q_{t-s}\right)(dz)ds
\label{eq:perdersecvera}
\\
&=\int_{0}^{t} \int_\calh
\left[\<\nabla\psi\left(
\bar f_1\left(s, (e^{sA}z)_0+(e^{tA}x)_0\right)\right),
s^{-{1/2}}
\bar{\bar f}\left(s, (e^{sA}z)_0+(e^{tA}x)_0\right) e^{(t-s)A}h\>_{\R^m} \right]
\caln\left(0,Q_{t-s}(dz) \right)ds
\nonumber
\end{align}
Now, from calculations similar to the ones performed in the first part we arrive at \red{
\begin{align}\label{eq:der-seconda-conv-Sigma}
\<\nabla^B\<\nabla\hat g(t,x),h\>_\calh,k\>_{\R^m}
&=\int_{0}^{t} \int_{\R^n}
\left[\<\nabla\psi\left(
\bar f_1\left(s, (e^{sA}z)_0+(e^{tA}x)_0\right)\right),
s^{-{1/2}}
\bar{\bar f}\left(s, z_1+(e^{tA}x)_0\right) e^{(t-s)A}h\>_{\R^m}
\right.
\nonumber
\\[2mm]
&
\qquad\qquad\left.
\<(Q^1_{t-s})^{-{1/2}}
\left( e^{tA} Bk\right)_0, (Q^1_{t-s})^{-{1/2}}z_1\>_{\R^n}\right]\caln\left(0,Q^1_{t-s}\right)(dz_1) ds.
\end{align}}
Since $\nabla \psi$ is bounded (as it satisfies \myref{eq:Hlip}) then (\ref{stimaiterata-secondag}) easily follows by the definition of $\bar{\bar f}$ and from \myref{eq:betaeuler}.
Second order differentiability and estimate (\ref{stimaiterata-secondagbis}), when $\sigma$ is onto, follow in the same way.
\qed



\section{Regular solutions of the HJB equation}
\label{sec-HJB}

We show first, in Subsection \ref{SS:EXUNMILD}, that the HJB equation (\ref{HJBformale}) admits a unique mild solution $v$ which is $B$-Fr\'echet
differentiable.
Then (Subsection \ref{SS:SECONDDERIVATIVE}) we prove a further regularity result
whose proof is more complicated than the previous one
and that will be useful to solve completely the control problem in the forthcoming companion paper \cite{FGFM-II}.

\subsection{Existence and uniqueness of mild solutions}
\label{SS:EXUNMILD}

We start
showing how to rewrite (\ref{HJBformale1}) in its
integral (or ``mild'') form as anticipated in the introduction, formula \myref{solmildHJB}.
Denoting by $\call$ the generator of the Ornstein-Uhlenbeck semigroup
$R_{t}$, we know that, for all $f\in C^{2}_b(\calh)$ such that $\nabla f \in D(A^*)$
(see e.g. \cite{CeGo} Section 5 or also \cite{DaPratoZabczyk95} Theorem 2.7):

\begin{equation}\label{eq:ell}
 \call[f](x)=\frac{1}{2} Tr \;GG^*\; \nabla^2f(x)
+ \< x,A^*\nabla f(x)\>.
\end{equation}
The HJB equation (\ref{HJBformale}) can then be formally rewritten as \red{
\begin{equation}\label{HJBformale}
  \left\{\begin{array}{l}\dis
-\frac{\partial v(t,x)}{\partial t}=\call [v(t,\cdot)](x) +\ell_0(t)+
H_{min} (\nabla^B v(t,x)),\qquad t\in [0,T],\,
x\in \calh,\\
\\
\dis v(T,x)=\bar\phi(x_0).
\end{array}\right.
\end{equation}
By applying formally the variation of
constants formula we then have
\begin{equation}
v(t,x) =R_{T-t}[\phi](x)+\int_t^T \left[R_{s-t}\left[
H_{min}(\nabla^B v(s,\cdot))\right](x)+\ell_0(s)\right]\; ds,\qquad t\in [0,T],\
x\in \calh,\label{solmildHJB}
\end{equation}}
We use this formula to give the notion of mild
solution for the HJB equation (\ref{HJBformale}).

\begin{definition}\label{defsolmildHJB}
We say that a
function $v:[0,T]\times \calh\rightarrow\mathbb{R}$ is a mild
solution of the HJB equation (\ref{HJBformale}) if the following
are satisfied:
\begin{enumerate}

\item $v(T-\cdot, \cdot)\in C^{0,1,B}_{{1/2}}\left([0,T]\times \calh\right)$;

\item  equality (\ref{solmildHJB}) holds on $[0,T]\times \calh$.
\end{enumerate}
\end{definition}

\begin{remark}\label{rm:crescitapoli-HJB}
Since $C^{0,1,B}_{{1/2}}\left([0,T]\times \calh\right)
\subset C_b([0,T]\times \calh)$ (see Definition \ref{df4:Gspaces})
the above Definition \ref{defsolmildHJB} requires, among other properties, that a mild solution is continuous and bounded up to $T$.
This constrains the assumptions on the data, e.g. it implies that the final datum $\phi$ must be continuous and bounded.
As recalled in Remark \ref{rm:crescitapoli-gen}-(i) and (ii)
we may change this requirement in the above definition asking only polynomial growth in $x$
and/or measurability of $\phi$. Most of our main results will remain true
with straightforward modifications.
\end{remark}

Since the transition semigroup $R_t$ is not even strongly Feller
we cannot study the existence and uniqueness of a mild solution of equation (\ref{HJBformale}) as it is done e.g. in \cite{G1}.
We then use the partial smoothing property studied in Sections \ref{section-smoothOU} and
\ref{section-smooth-conv}.
Due to Lemma \ref{lemma reg-convpercontr} the right space where
to seek a mild solution seems to be $\Sigma^1_{T,{1/2}}$;
indeed our existence and uniqueness result will be proved by a fixed point argument in such space.

\begin{theorem}\label{esistenzaHJB}
Let Hypotheses \ref{ipotesibasic} and \ref{ipotesicostoconcreto} hold and
let (\ref{eq:hpdebreg}) or (\ref{eq:hpdebregbis}) hold.
Then the HJB equation (\ref{HJBformale})
admits a mild solution $v$ according to Definition \ref{defsolmildHJB}.
Moreover $v$
is unique among the functions $w$ such that $w(T-\cdot,\cdot)\in\Sigma_{T,1/2}$ and it satisfies, for suitable $C_T>0$, the estimate
\begin{equation}\label{eq:stimavmainteo}
\Vert v(T-\cdot,\cdot)\Vert_{C^{0,1,B}_{{1/2}}}\le C_T\left(\Vert\bar\phi \Vert_\infty
+\Vert\bar\ell_0 \Vert_\infty \right).
\end{equation}
Finally if the initial datum $\phi$ is also continuously $B$-Fr\'echet
(or Fr\'echet) differentiable,
then $v \in C^{0,1,B}_{b}([0,T]\times \calh)$ and, for suitable $C_T>0$,
\begin{equation}\label{eq:stimavmainteobis}
\Vert v\Vert_{C^{0,1,B}_{b}}\le C_T\left(\Vert\phi \Vert_\infty
+\Vert\nabla^B\phi \Vert_\infty+\Vert\ell_0 \Vert_\infty \right)
\end{equation}
(substituting $\nabla^B\phi$ with $\nabla\phi$ if $\phi$ is Fr\'echet differentiable).
\end{theorem}
\dim
We use a fixed point argument in $\Sigma^1_{T,{1/2}}$. To this aim,
first we rewrite (\ref{solmildHJB}) in a forward way. Namely
if $v$ satisfies \myref{solmildHJB} then, setting $w(t,x):=v(T-t,x)$ for any
$(t,x)\in[0,T]\times \calh$, we get that $w$ satisfies 
\begin{equation}
 w(t,x) =R_{t}[\phi](x)+\int_0^t\left[ R_{t-s}[H_{min}(\nabla^B w(s,\cdot))](x)+\red{\ell_0(s)}\right]\; ds,\qquad t\in [0,T],\,x\in \calh,\label{solmildHJB-forward}
\end{equation}
which is the mild form of the forward HJB equation
\begin{equation}\label{HJBformaleforward}
  \left\{\begin{array}{l}\dis
\frac{\partial w(t,x)}{\partial t}=\call [w(t,\cdot)](x) +\red{\ell_0(t)}+
H_{min} (\nabla^B w(t,x)),\qquad t\in [0,T],\,
x\in \calh,\\
\\
\dis w(0,x)=\phi(x).
\end{array}\right.
\end{equation}
Define the map $\calc$ on $\Sigma^1_{T,{1/2}}$  by setting, for $g\in \Sigma^1_{T,{1/2}}$,
\begin{equation}\label{mappaC}
 \calc(g)(t,x):=R_{t}[\phi](x)+\int_0^t \red{\left[R_{t-s}[
H_{min}(
\nabla^B g(s,\cdot))
](x)+\ell_0(s)\right]}\; ds,\qquad t\in [0,T],\,
x\in \calh.
\end{equation}
By Proposition \ref{lemmaderhpdeb}
and Lemma \ref{lemma reg-convpercontr}-(i) we deduce that $\calc$
is well defined in $\Sigma^1_{T,{1/2}}$ and takes its values in $\Sigma^1_{T,{1/2}}$. Since in Lemma \ref{lemma:Sigma}
we have proved that $\Sigma^1_{T,{1/2}}$ is a closed subspace of $C^{0,1,B}_{{1/2}}([0,T]\times\calh)$, once we have proved that $\calc$ is a contraction,
by the Contraction Mapping Principle there exists a unique (in $\Sigma^1_{T,{1/2}}$) fixed point of the map $\calc$, which gives a mild solution of (\ref{HJBformale}).

\noindent Let $g_1,g_2 \in \Sigma^1_{T,{1/2}}$. We evaluate
$\Vert \calc(g_1)-\calc (g_2)\Vert_{\Sigma_{T,{1/2}}}=\Vert \calc(g_1)-\calc (g_2)\Vert_{C^{0,1,B}_{{1/2}}}$. First of all, arguing as in the proof of
Lemma \ref{lemma reg-convpercontr} we have, for every $(t,x)\in [0,T]\times \calh$,
\begin{align*}
  \vert \calc (g_1)(t,x)- \calc(g_2)(t,x) \vert& =\left\vert \int_0^t R_{t-s}\left[H_{min}\left(\nabla^B g_1(s,\cdot)\right)
 -H_{min}\left(\nabla^B g_2(s,\cdot)\right)\right](x)ds\right\vert\\
 &\le \int_0^t s^{-{1/2}} L \sup_{y \in H}\vert s^{{1/2}}\nabla^B (g_1-g_2)(s,y)\vert ds
  \leq 2Lt^{{1/2}}\Vert g_1-g_2 \Vert_{C^{0,1,B}_{{1/2}}}.
\end{align*}
Similarly, arguing exactly as in the proof of (\ref{stimaiterata-primag}), we get
\begin{align*}
t^{{1/2}}\vert \nabla^B\calc (g_1)(t,x) &- \nabla^B\calc(g_2)(t,x) \vert =
t^{{1/2}}\left\vert \nabla^B\int_0^t  R_{t-s}\left[H_{min}
\left(\nabla^B g_1(s,\cdot)\right)-H_{min}\left(\nabla^B
g_2(s,\cdot)\right)\right](x)ds\right\vert\\
& \leq t^{{1/2}} L \Vert g_1-g_2 \Vert_{C^{0,1,B}_{{1/2}}} \int_0^t (t-s)^{-{1/2}}
s^{-{1/2}} ds
\le t^{{1/2}}L \beta\left({1/2} , {1/2}\right)\Vert g_1-g_2 \Vert_{C^{0,1,B}_{{1/2}}}.
\end{align*}
Hence, if $T$ is sufficiently small, we get
\begin{equation}\label{stima-contr}
 \left\Vert \calc (g_1)-\calc(g_2)\right\Vert _{C^{0,1,B}_{{1/2}}
  }\leq C
\left\Vert g_1-g_2\right\Vert _{C^{0,1,B}_{{1/2}}
}
\end{equation}
with $C<1$. So the map $\calc$ is a contraction in $\Sigma^1_{T,{1/2}}$
and, if we denote by $w$ its unique fixed point, then $v:=w(T-\cdot,\cdot)$
turns out to be a mild solution of the HJB equation (\ref{HJBformale}),
according to Definition \ref{solmildHJB}.

Since the constant $L$ is independent of $t$, the case of generic $T>0$ follows
by dividing the interval $[0,T]$
into a finite number of subintervals of length $\delta$ sufficiently small, or equivalently, as done in \cite{Mas},
by taking an equivalent norm with an adequate exponential weight, such as
\[
 \left\Vert f\right\Vert _{\eta,C^{0,1,B}_{{1/2}}
 }=\sup_{(t,x)\in[0,T]\times \calh}
\vert e^{\eta t}f(t,x)\vert+
\sup_{(t,x)\in (0,T]\times \calh}  e^{\eta t}t^{{1/2}}
\left\Vert \nabla^B f\left(  t,x\right)  \right\Vert _{(\R^m)^{\ast}},
\]


The estimate (\ref{eq:stimavmainteo}) follows from Proposition \ref{cor-der} and Lemma \ref{lemma reg-convpercontr}.

Finally the proof of the last statement follows observing that, if
$\phi$ is continuously $B$-Fr\'echet (or Fr\'echet) differentiable,
then $R_t[\phi]$ is continuously $B$-Fr\'echet differentiable
with $\nabla^B R_t[\phi]$ bounded in $[0,T]\times \calh$, see lemma \ref{lemma-reg-R_t},
formula (\ref{eq:derivateprimelisce}).
This allows to perform the fixed point, exactly as done in the first part of the proof,
in $C^{0,1,B}_b([0,T]\times \calh)$ and to prove estimate \myref{eq:stimavmainteobis}.
\qed

\begin{corollary}\label{diffle-corollario}
Let Hypotheses \ref{ipotesibasic} and \ref{ipotesicostoconcreto} hold and let $\sigma$ be onto.
Then the mild solution of equation (\ref{HJBformale})
found in the previous theorem
is also Fr\'echet differentiable,
and the following estimate holds true
\begin{equation}\label{stimadiffle-solHJB}
 \left\Vert v(T-\cdot,\cdot)\right\Vert _{C^{0,1}_{{1/2}}
  }\leq C_T\left(\Vert\bar\phi \Vert_\infty
+\Vert\bar\ell_0 \Vert_\infty \right)
\end{equation}
for a suitable $C_T>0$.
\end{corollary}
\dim
Let $v$ be the mild solution of equation (\ref{HJBformale}), and $\forall
t\in[0,T],\,x\in\calh$ define $w(t,x):=v(T-t,x)$, so that $w$ satisfies
(\ref{solmildHJB-forward}), so that by applying
the last statement of Lemma \ref{lemma reg-convpercontr}
it is immediate to see that $w\in C^{0,1}_{{1/2}}([0,T]\times \calh)$.
By differentiating (\ref{solmildHJB-forward}) we get
\begin{equation*}
 \nabla w(t,x) =\nabla R_{t}[\phi](x)+\nabla\int_0^t \left[R_{t-s}[H_{min}(\nabla^B w(s,\cdot)](x)+\red{\ell_0(s)}\right]\; ds,\qquad t\in [0,T],\,x\in \calh,
\end{equation*}
By Lemma \ref{lemma-reg-R_t}, the above recalled variation of Lemma
\ref{lemma reg-convpercontr} and estimate (\ref{stimaiterata-primagbis}), we get that
\begin{equation*}
 \vert\nabla w(t,x) \vert\leq C t^{-{1/2}}\Vert\phi\Vert_\infty+Ct^{{1/2}}\left(1+\Vert
 w\Vert_{C^{0,1,B}_{{1/2}}}+\Vert\ell_0\Vert_\infty\right),
\qquad t\in [0,T],\
x\in \calh,
\end{equation*}
which gives the claim using the estimate for
$\Vert w\Vert_{C^{0,1,B}_{{1/2}}}$ given in (\ref{eq:stimavmainteo}).
\qed

\subsection{Second derivative of mild solutions}
\label{SS:SECONDDERIVATIVE}

The further regularity result we are going to prove
is interesting in itself, but is also crucial
to solve the control problem, as will be seen in
the forthcoming companion paper \cite{FGFM-II}.
A similar result can be found in \cite{G1}, Section 4.2. Here we use the same line of proof
but we need to argue in a different way to get the apriori estimates.

\begin{theorem}\label{lemma-stimev}
Let Hypotheses \ref{ipotesibasic}, \ref{ipotesicostoconcreto} and
\ref{ipotesicostoconcretobis} hold.
Let also (\ref{eq:hpdebreg}) or (\ref{eq:hpdebregbis}) hold.
Let $v$ be the mild solution of the HJB equation (\ref{HJBformale}) as from
Theorem \ref{esistenzaHJB}. Then we have the following.
\begin{itemize}
  \item[(i)]
If $\phi$ is continuously differentiable then we have $v \in \Sigma^2_{T,{1/2}}$,
hence the second order derivatives $\nabla^B\nabla v$ and $\nabla\nabla^B v$
exist and are equal on $[0,T)\times \calh$.
Moreover there exists a constant $C>0$ such that, on $[0,T)\times \calh$,
\red{\begin{equation}\label{stimanablav}
\vert \nabla v(t,x)\vert \leq C \Vert \nabla\bar\phi\Vert_\infty,
\end{equation}
\begin{equation}\label{stimanablav^2}
\vert \nabla^B\nabla v(t,x)\vert= \vert \nabla\nabla^B v(t,x)\vert
\leq C (T-t)^{-{1/2}}
\Vert \nabla\bar\phi\Vert_\infty.
\end{equation}
Finally, if $\sigma$ is onto, then also $\nabla^2 v$ exists and is continuous
on on $[0,T)\times \calh$ and, on such set, for suitable $C>0$,
\begin{equation}\label{stimanablav^2bis}
\vert \nabla^2 v(t,x)\vert\leq C (T-t)^{-{1/2}}\Vert \nabla\bar\phi\Vert_\infty.
\end{equation}}

  \item[(ii)]
If $\phi$ is only continuous then the function
$(t,x)\mapsto (T-t)^{1/2}v(t,x)$ belongs to $\Sigma^2_{T,{1/2}}$.
Moreover there exists a constant $C>0$ such that, on $[0,T)\times \calh$, \red{
\begin{equation}\label{stimanablavreg}
\vert \nabla v(t,x)\vert \leq C(T-t)^{-1/2} \Vert\bar\phi\Vert_\infty,
\end{equation}
\begin{equation}\label{stimanablav^2reg}
\vert \nabla^B\nabla v(t,x)\vert= \vert \nabla\nabla^B v(t,x)\vert\leq C (T-t)^{-1}\Vert \bar\phi\Vert_\infty.
\end{equation}
Finally, if $\sigma$ is onto, then also $\nabla^2 v$ exists and is continuousin $[0,T)\times \calh$ and on such set, for suitable $C>0$,
\begin{equation}\label{stimanablav^2bisreg}
\vert \nabla^2 v(t,x)\vert\leq C(T-t)^{-1}\Vert \nabla\bar\phi\Vert_\infty.
\end{equation}}
\end{itemize}
\end{theorem}	
\dim
We start proving (i) by applying the Contraction Mapping Theorem in a closed ball $B_T(0,R)$
($R$ to be chosen later) of the space $ \Sigma^2_{T,{1/2}}$.
By Proposition \ref{lemma-reg-R_t} and Lemma \ref{lemma reg-convpercontr} we deduce that
the map $\calc$ defined in \myref{mappaC}
brings $\Sigma^2_{T,{1/2}}$ into $\Sigma^2_{T,{1/2}}$.
Moreover, for every $g \in \Sigma^2_{T,{1/2}}$, we get, first using \myref{eq:Hlip},
\begin{align*}
&  \vert \calc (g)(t,x) \vert\le
\left\vert
R_{t}[\phi](x)
\right\vert
+
\left\vert \int_0^t \left[R_{t-s}\left[H_{min}\left(\nabla^B g(s,\cdot)\right)
\right](x)+\red{ \ell_0(s)}\right]ds
\right\vert
\leq
\Vert \phi\Vert_\infty
+t L\left(1+\Vert g \Vert_{C^{0,1,B}_{0}}\right) + t\|\ell_0\|_\infty;
\end{align*}
second by \myref{eq:derivateprimelisce}, \myref{eq:perdersecvera}, \myref{eq:Hlip}
(calling $M:=\sup_{[0,T]}\|e^{tA}\|$)\red{ 
\begin{align}
\nonumber
  \vert\nabla \calc (g)(t,x) \vert &\le
\left\vert
\nabla R_{t}[\phi](x)
\right\vert
+
\left\vert \nabla\int_0^t \left[R_{t-s}\left[H_{min}\left(\nabla^B g(s,\cdot)\right)\right](x)+\ell_0(s)\right]ds
\right\vert\\[2mm]
&\le M \Vert \nabla\phi\Vert_\infty+M
\int_0^t \left[\left\Vert\nabla H_{min}\left(\nabla^B g(s,\cdot)\right)
\nabla\nabla^B g(s,\cdot)\right\Vert_\infty \right]ds
\nonumber
\\
& \leq M\left[ \Vert \nabla\phi\Vert_\infty 
 + L t^{{1/2}}\Vert g \Vert_{C^{0,2,B}_{{1/2}}}\right];
\label{eq:stimadersecnew1}
\end{align}
third by \myref{eq:stimaDRT2-gen} (with
(\ref{eq:hpdebreg}) or (\ref{eq:hpdebregbis})), and \myref{stimaiterata-secondag}
\begin{align}
\nonumber
&t^{{1/2}}\vert \nabla^B\nabla\calc (g)(t,x)\vert
\le
t^{{1/2}}
\left\vert
\nabla^B\nabla R_{t}[\phi](x)
\right\vert
+
t^{{1/2}} \left\vert \nabla^B\nabla
\int_0^t  R_{t-s}\left[H_{min}\left(\nabla^B g(s,\cdot)\right)
\right](x)ds\right\vert
\\[2mm]
& \leq C \Vert \nabla\phi\Vert_\infty
+ Ct^{{1/2}} \| g\|_{C^{0,2,B}_{1/2}}
\leq C\left[ \Vert \nabla\phi\Vert_\infty
+t^{{1/2}} \Vert g \Vert_{C^{0,2,B}_{{1/2}}}\right],
\label{eq:stimadersecnew2}
\end{align}
with the constant $C$ (that may change from line to line) given by the quoted estimates.
Hence, for $g\in B_T(0,R)$, we get, for given $C_1>0$,
\begin{align}
\nonumber
\Vert \calc (g) \Vert_{C^{0,2,B}_{{1/2}}}
&\le
C_1\left[\Vert\phi\Vert_{C^1_b} +T \right]
+TL \Vert g \Vert_{C^{0,1,B}_{0}}+(ML+C)T^{1/2}\Vert g \Vert_{C^{0,2,B}_{{1/2}}}
\\
\label{eq:dersecstima1}
&\le C_1\left[\Vert\phi\Vert_{C^1_b} +T\right]+\rho(T)R
\end{align}}
where we define
\begin{equation}\label{eq:rhodef}
\rho(T):= TL+(ML+C)T^{1/2}.
\end{equation}
Now take $g_1,g_2 \in \Sigma^2_{T,{1/2}}$.
Arguing as in the above estimates we have, for every $(t,x)\in [0,T]\times \calh$,
\begin{align*}
\vert \calc (g_1)(t,x)- \calc(g_2)(t,x) \vert&=\left\vert \int_0^t R_{t-s}\left[H_{min}\left(\nabla^B g_1(s,\cdot)\right)
-H_{min}\left(\nabla^B g_2(s,\cdot)\right)\right](x)ds\right\vert\\[2mm]
&\leq t L\Vert g_1-g_2 \Vert_{C^{0,1,B}_{0}}
\end{align*}
\begin{align*}
&\vert \nabla\calc (g_1)(t,x)- \nabla\calc(g_2)(t,x) \vert =
\left\vert \nabla\int_0^t  R_{t-s}\left[H_{min}\left(\nabla^B g_1(s,\cdot)\right)
-H_{min}\left(\nabla^B g_2(s,\cdot)\right)\right](x)ds
\right\vert
\\
&\le M\int_0^t \left\Vert\nabla H_{min}\left(\nabla^B g_1(s,\cdot)\right)
\nabla\nabla^B g_1(s,\cdot)-
\nabla H_{min}\left(\nabla^B g_2(s,\cdot)\right)\nabla\nabla^B g_2(s,\cdot)\right\Vert_\infty
ds
\\
& \leq  2ML t^{{1/2}}\left[\Vert g_1-g_2 \Vert_{C^{0,1,B}_{0}}
\Vert g_1 \Vert_{C^{0,2,B}_{{1/2}}}+\Vert g_1-g_2 \Vert_{C^{0,2,B}_{{1/2}}} \right]
\end{align*}
and, using \myref{eq:der-seconda-conv-Sigma},
\begin{align*}
&t^{{1/2}}\vert \nabla^B\nabla\calc (g_1)(t,x)- \nabla^B\nabla\calc(g_2)(t,x) \vert =
t^{{1/2}} \left\vert \nabla^B\nabla\int_0^t  R_{t-s}\left[H_{min}\left(\nabla^B g_1(s,\cdot)\right)
-H_{min}\left(\nabla^B g_2(s,\cdot)\right)\right](x)ds\right\vert
\\[2mm]
& \leq  t^{{1/2}}M L \beta\left({1/2} , {1/2}\right)
\left[\Vert g_1-g_2 \Vert_{C^{0,1,B}_{0}}
\Vert g_1 \Vert_{C^{0,2,B}_{{1/2}}}+\Vert g_1-g_2 \Vert_{C^{0,2,B}_{{1/2}}} \right].
%
\end{align*}
Hence, for $g_1,g_2\in B_T(0,R)$, we have, recalling \myref{eq:rhodef} and the way $C$ is found
in \myref{eq:stimadersecnew2},
\begin{align}
\nonumber
\left\Vert \calc (g_1)-\calc(g_2)\right\Vert _{C^{0,2,B}_{{1/2}}([0,T]\times \calh)  }
&\leq L\left(T+ MT^{{1/2}}\left(2+\beta\left({1/2} , {1/2}\right)\right)
(1+R) \right)
\left\Vert g_1-g_2\right\Vert _{C^{0,2,B}_{{1/2}}  }
\\[2mm]
&\le \rho(T)(1+R)\left\Vert g_1-g_2\right\Vert _{C^{0,2,B}_{{1/2}}  },
\label{eq:dersecstima2}
\end{align}
Now, by \myref{eq:dersecstima1} and \myref{eq:dersecstima2}, choosing any \red{$R>C_1\left[\Vert\phi\Vert_{C^1_b} +T\right]$} we can find $T_0$ sufficiently small so that $\rho(T_0)<1/2$ and so, thanks to \myref{eq:dersecstima1} and \myref{eq:dersecstima2}, $\calc$ is a contraction in $B_{T_0}(0,R)$.
Let then $w$ be the unique fixed point of $\calc$ in $B_{T_0}(0,R)$:
it must coincide with $v(T-\cdot,\cdot)$ for $t\in [0,T_0]$.
This procedure can be iterated arriving to cover the whole interval $[0,T]$
if we give an apriori estimate for the norm $\Vert w \Vert_{C^{0,2,B}_{{1/2}}}$.
By the last statement of Theorem \ref{esistenzaHJB} we already have an apriori estimate for
$\Vert w \Vert_\infty+\Vert \nabla^B w \Vert_\infty$.
To get the estimate for $\nabla w$ and $\nabla^B\nabla w $
we use \myref{eq:stimadersecnew1} and \myref{eq:stimadersecnew2}
where we put $w$ in place of $\calc g$ and $g$.
From the first line of \myref{eq:stimadersecnew2} and \ref{eq:der-seconda-conv-Sigma}
we get \red{
$$
\Vert \nabla^B \nabla w(t,\cdot)\Vert_\infty
\leq C\Vert \nabla\phi\Vert_\infty
+L\int_0^t (t-s)^{-{1/2}}
\Vert \nabla^B \nabla w(s,\cdot)\Vert_\infty ds
$$}
which, thanks to the Gronwall Lemma  (see \cite{Henry}, Subsection 1.2.1, p.6)
give the apriori estimate for $\nabla^B \nabla w$. Then from the second line of
\myref{eq:stimadersecnew1} we get \red{
$$
\Vert \nabla w(t,\cdot)\Vert_\infty
\le  M\Vert \nabla\phi\Vert_\infty 
+M L\int_0^t \Vert \nabla^B \nabla w(s,\cdot)\Vert_\infty ds
$$}
which gives the apriori estimate for $\nabla w$ using the previous one for $\nabla^B \nabla w$.
Estimate \myref{stimanablav^2bis} follows by repeating the same arguments above
but replacing $\nabla^B \nabla $ with $\nabla^2$.

We now prove (ii).
Let $v$ be the mild solution of
(\ref{HJBformale}) and, for all $\eps \in ]0, T[$, $x \in \calh$, call $\phi^\eps (x) =
v(T-\eps,x)$. Then $v$ is the unique mild solution, on $[0,T-\eps]\times \calh$,
of the equation (for $t \in [0,T]$, $x \in \calh$) 
\begin{equation}
\label{eq:HJBmildeps}
v(t,x) =   R_{T-\eps-t}  \phi^\eps (x)+ \int_{t}^{T-\eps}\left[
R_{s-t}\left[H_{min}(\nabla^B v(s,\cdot ))\right](x)+\red{+\ell_0(s)}\right] ds
\end{equation}
This fact can be easily seen by applying the semigroup property of $R_t$
(see e.g. \cite{G1} Lemma 4.10 for a completely similar result).

Now, by Theorem \ref{esistenzaHJB}, $\phi^\eps$ is continuously differentiable,
so we can apply part (i) of this theorem to \myref{eq:HJBmildeps}
getting the required $C^2$ regularity.
Estimates
\myref{stimanablavreg}-\myref{stimanablav^2reg}-\myref{stimanablav^2bisreg}
follows using estimates
\myref{stimanablav}-\myref{stimanablav^2}-\myref{stimanablav^2bis}
with $\phi^\eps$ in place of $\phi$ and then using the arbitrariness of $\eps$ and applying
\myref{eq:stimavmainteo} to estimate $\phi^\eps$ in term of $\phi$.
\qed

We emphasize that Theorem \ref{lemma-stimev} provides more regularity for the mild solution of the HJB equation (\ref{HJBformale}) when the data are more regular.
In particular, when $\sigma$ is onto, the classical first and second derivatives of the solution exist and are continuous. In such case it seems also possible to prove, on the line of what is done e.g. in Theorem 7.5.1 of \cite{DP3}, that the mild solution is also a strict solution of the HJB equation; i.e., roughly speaking, that it admits time derivative and solves \myref{HJBformale1} pointwise in $[0,T]\times D(A)$.
We do not do it here as our main goal is solve the control problem proving a verification theorem and the existence of optimal feedbacks, as we do in  the companion paper \cite{FGFM-II}, Section 5.
To accomplish this task is indeed enough to prove the weaker fact that the mild solution can be approximated by a sequence of strict solutions of suitable approximating equations, as we do in \cite{FGFM-II}, Section 4, see in particular Remark 4.5 there.

On the other hand it does not seem possible, as explained in the companion paper \cite{FGFM-II}, Remark 4.4, to prove the existence of an approximating sequence of classical solutions,
where by classical solution we mean a regular solution which satisfies \myref{HJBformale}
pointwise in $[0,T]\times \calh$. This suggest that, even when the mild solution is also a strict solution, it may not be a classical solution: the point is
that we cannot guarantee that a smooth solution $v$ of the HJB equation
(\ref{HJBformale}) is such that $\nabla v\in D(A^*)$.

\end{document}